\documentclass[12pt]{article}
\usepackage{latexsym, amsthm, amsmath, amssymb}
\newcommand{\duk}{\noindent {\bf Proof. }}
\newcommand{\kduk}{\hfill $\Box$\bigskip}
\newcommand{\R}{\mathbb{R}}
\newcommand{\N}{\mathbb{N}}
\newcommand{\No}{\mathbb{N}_0}
\newcommand{\Z}{\mathbb{Z}}
\newcommand{\Q}{\mathbb{Q}}
\newcommand{\C}{\mathbb{C}}
\newcommand{\bsn}{\bigskip\noindent}
\newcommand{\sus}{\subset}
\def\ds{\dots}
\def\av{\mathrm{Av}}
\def\ep{\varepsilon}
                                          
\newtheorem{lemma}{Lemma}[section]

\newtheorem{theorem}[lemma]{Theorem}
\newtheorem{proposition}[lemma]{Proposition}

\newtheorem{problem}[lemma]{Problem}

\title{Overview of some general results in combinatorial enumeration}
\author{Martin Klazar\thanks{Charles University, Faculty of Mathematics and 
Physics, Department of Applied Mathematics (KAM) and Institute for Theoretical Computer 
Science (ITI), Malostransk\'e n\'am. 25, 
Praha 11800, Czech Republic. ITI is supported by the project 1M0021620808 of the Czech 
Ministry of Education. Email: klazar at kam.mff.cuni.cz.}}

\begin{document}\maketitle

\begin{abstract}
This survey article is devoted to general results in combinatorial enumeration. The 
first part surveys results on growth of hereditary properties of combinatorial structures.
These include permutations, ordered and unordered graphs and hypergraphs, relational 
structures, and others. The second part advertises five topics in general 
enumeration: 1. counting lattice points in lattice polytopes,   
2. growth of context-free languages, 3. holonomicity 
(i.e., P-recursiveness) of numbers of labeled regular graphs, 4. frequent occurrence of the 
asymptotics $cn^{-3/2}r^n$ and 5. 
ultimate modular periodicity of numbers of MSOL-definable structures.
\end{abstract}

\section{Introduction}

We survey some general results in combinatorial enumeration. 
A {\em problem} in enumeration is (associated with) an infinite sequence 
$P=(S_1,S_2,\dots)$ of finite sets $S_i$. Its {\em counting function} $f_P$ is given by 
$f_P(n)=|S_n|$, the cardinality of the set $S_n$. We are 
interested in results of the following kind on {\em general} classes of problems and their counting functions. 

\bsn 
{\bf Scheme of general results in combinatorial enumeration.} {\em The counting 
function $f_P$ of every problem $P$ in the class ${\cal C}$ belongs to the class of functions ${\cal F}$. Formally, $\{f_P\;|\;P\in{\cal C}\}\sus{\cal F}$.}

\bsn
The larger ${\cal C}$ is, and the more specific the functions in ${\cal F}$ are, 
the stronger the result. The present overview is a collection of many examples of this scheme. 

One can distinguish general results of two types. In {\em exact results},  
${\cal F}$ is a class of explicitly defined functions, for example polynomials or functions defined by recurrence relations of certain type or 
functions computable in polynomial time. 
In {\em asymptotic results}, ${\cal F}$ consists of functions defined by asymptotic equivalences or asymptotic  inequalities, for example functions growing at most exponentially or functions asymptotic to $n^{(1-1/k)n+o(n)}$ as $n\to\infty$, with the
constant $k\ge 2$ being an integer.

The sets $S_n$ in $P$ usually constitute sections of a fixed infinite
set. Generally speaking, we take an infinite universe $U$ of combinatorial structures and 
introduce problems and classes of problems as subsets of $U$ and 
families of subsets of $U$, by means of {\em size functions} $s:\;U\to\No=\{0,1,2,\ds\}$ and/or 
(mostly binary) relations between structures in $U$. More specifically,  
we will mention many results falling within the framework of growth of downsets in partially order sets, or posets.

{\bf Downsets in posets of combinatorial structures.} We consider a 
nonstrict partial ordering $(U,\prec)$, where $\prec$ is a containment or a 
substructure relation on a set $U$ of combinatorial structures, and a size function 
$s:\;U\to\No$. Problems $P$ are {\em downsets} in $(U,\prec)$, meaning that $P\sus U$ 
and $A\prec B\in P$ implies $A\in P$, and the counting function of $P$ is
$$
f_P(n)=\#\{A\in P\;|\;s(A)=n\}.
$$  
(More formally, the problem is the sequence of sections $(P\cap U_1,P\cap U_2,\dots)$ where $U_n=\{A\in U\;|\;s(A)=n\}$.) 
Downsets are exactly the sets of the form
$$
\av(F):=\{A\in U\;|\;A\not\succ B\mbox{ for every $B$ in $F$}\},\ F\sus U.
$$
There is a one-to-one
correspondence $P\mapsto F=\min(U\backslash P)$ and $F\mapsto P=\av(F)$ between the 
family of downsets $P$ and the family of {\em antichains} $F$, which are sets of mutually incomparable structures under $\prec$. We call the antichain $F=\min(U\backslash P)$ corresponding 
to a downset $P$ the {\em base of $P$}.

We illustrate the scheme by three examples, all for downsets in posets. 

\subsection{Three examples}

{\bf Example 1. Downsets of partitions.} $U$ is the family of partitions of 
$[n]=\{1,2,\dots,n\}$ for $n$ ranging in $\N$, so $U$ consists of finite sets 
$S=\{B_1, B_2,\ds,B_k\}$ of disjoint and nonempty finite subsets $B_i$ of $\N$, 
called {\em blocks}, whose union  
$B_1\cup B_2\cup\ds\cup B_k=[n]$ for some $n$ in $\N$. 
Two natural size functions on $U$ are order and size, where the {\em order}, $\|S\|$, of $S$
is the cardinality, $n$, of the underlying set and the
{\em size}, $|S|$, of $S$ is the number, $k$, of blocks. The formula for the number of 
partitions of $[n]$ with $k$ blocks
$$
S(n,k):=\#\{S\in U\;|\;\|S\|=n, |S|=k\}=\sum_{i=0}^k\frac{(-1)^i(k-i)^n}{i!(k-i)!}
$$
is a classical result (see \cite{stan1}); $S(n,k)$ are called {\em Stirling numbers}. It is 
already a simple example of the above scheme but we shall go further. 

For fixed $k$, the function $S(n,k)$ is a linear combination with 
rational coefficients of the exponentials $1^n,2^n,\ds,k^n$. So is the sum 
$S(n,1)+S(n,2)+\ds+S(n,k)$ counting partitions with order $n$ and size at most $k$.
We denote the set of such partitions $\{S\in U\;|\;|S|\le k\}$ as $U_{\le k}$. 
Consider the poset $(U,\prec)$ with $S\prec T$ meaning that there
is an {\em increasing} injection $f:\;\bigcup S\to\bigcup T$ such that every two elements 
$x,y$ in $\bigcup S$ lie in the same block of $S$ if and only if $f(x),f(y)$ lie in 
the same block of $T$. In other words, $S\prec T$ means that $\bigcup T$ has a subset 
$X$ of size $\|S\|$ such that $T$ induces on $X$ a partition order-isomorphic to $S$.
Note that $U_{\le k}$ is a downset in $(U,\prec)$. We know that the counting 
function of $U_{\le k}$ with respect to order $n$ equals $a_11^n+\ds+a_kk^n$ with $a_i$ in 
$\Q$. What are the counting functions of other downsets? If the size is bounded, as for 
$U_{\le k}$, they have similar form as shown in the next theorem, proved by Klazar \cite{klaz_stir}. It is our first example of an exact general enumerative result. 

\begin{theorem}[Klazar]\label{setpartstir}
If $P$ is a downset in the poset of partitions such that   
$\max_{S\in P}|S|=k$, then there exist a natural number $n_0$ and polynomials $p_1(x)$, $p_2(x)$, \ldots, $p_k(x)$ with rational coefficients such that for every $n>n_0$,
$$
f_P(n)=\#\{S\in P\;|\;\|S\|=n\}=p_1(n)1^n+p_2(n)2^n+\ds+p_k(n)k^n.
$$
\end{theorem}

\bsn
If $\max_{S\in P}|S|=+\infty$, the situation is 
much more intricate and we are far from having a complete description but the growths of 
$f_P(n)$ below $2^{n-1}$ have been determined (see Theorem~\ref{BBMorgr} and the 
following comments). We briefly mention three 
subexamples of downsets with unbounded size, none of which has $f_P(n)$ in the form 
of Theorem~\ref{setpartstir}. 
If $P$ consists of all partitions of $[n]$ into intervals of length at most $2$, then $f_P(n)=F_n$, the 
$n^{\mbox{\scriptsize{th}}}$ Fibonacci number, and $f_P(n)=b_1\alpha^n+b_2\beta^n$
where $\alpha=\frac{\sqrt{5}-1}{2}$, $\beta=\frac{\sqrt{5}+1}{2}$ and 
$b_1=\frac{\alpha}{\sqrt{5}},b_2=\frac{\beta}{\sqrt{5}}$. If $P$ is given
as $P=\av(\{C\})$ where $C=\{\{1,3\},\{2,4\}\}$ (the partitions in $P$ are so called 
{\em noncrossing partition}, see the survey of Simion \cite{simi}) then 
$f_P(n)=\frac{1}{n+1}{2n\choose n}$, the $n^{\mbox{\scriptsize{th}}}$ Catalan 
number which is asymptotically $cn^{-3/2}4^n$. Finally, if $P=U$, so $P$ consists of all 
partitions, then $f_P(n)=B_n$, the $n^{\mbox{\scriptsize{th}}}$ Bell number which grows superexponentially.

{\bf Example 2. Hereditary graph properties.}
$U$ is the universe of finite simple graphs $G=([n],E)$ with vertex
sets $[n]$, $n$ ranging over $\N$, and $\prec$ is the induced subgraph relation;  
$G_1=([n_1],E_1)\prec G_2=([n_2],E_2)$ means that there is an injection from $[n_1]$ to $[n_2]$ (not necessarily increasing) that sends edges to edges and nonedges to nonedges.
The size, $|G|$, of a graph $G$ is the number of vertices. Problems are downsets in 
$(U,\prec)$ and are called {\em hereditary graph properties}. The next theorem, proved 
by Balogh, Bollob\'as and Weinreich \cite{balo_boll_wein_spee}, 
describes counting functions of hereditary graph properties that grow no faster than exponentially.

\begin{theorem}[Balogh, Bollob\'as and Weinreich]\label{herprostir}
If $P$ is a hereditary graph property such that for some constant $c > 1$, $f_P(n)=\#\{G\in P\;|\;|G|=n\}<c^n$ for 
every $n$ in $\N$, then there exists a natural number $n_0$ and  polynomials $p_1(x)$, $p_2(x)$, \ldots, $p_k(x)$ with rational coefficients such that for every $n>n_0$,
$$
f_P(n)=p_1(n)1^n+p_2(n)2^n+\ds+p_k(n)k^n.
$$
\end{theorem}

\noindent
The case of superexponential growth of $f_P(n)$ is discussed below in 
Theorem~\ref{BBWherpro}. 

In both examples we have the same class of functions ${\cal F}$, linear combinations 
$p_1(n)1^n+p_2(n)2^n+\ds+p_k(n)k^n$ with $p_i\in\Q[x]$. It would be nice to find a common extension of Theorems~\ref{setpartstir} and \ref{herprostir}.
It would be also of interest to determine if the two classes of
functions realizable as counting functions in both theorems coincide and how they 
differ from $\Q[x,2^x,3^x,\ds]$.

{\bf Example 3. Downsets of words.} $U$ is the set of finite 
words over a finite alphabet $A$, so $U=\{u=a_1a_2\dots a_k\;|\;a_i\in A\}$. 
The size, $|u|$, of such a word is its length $k$. The subword relation
$u=a_1a_2\ds a_k\prec v=b_1b_2\ds b_l$ means that 
$b_{i+1}=a_1,b_{i+2}=a_2,\ds,b_{i+k}=a_k$ for some $i$. We associate with 
an infinite word $v=b_1b_2\dots$ over $A$ the set $P=P_v$ of all its finite subwords, 
thus $P_v=\{b_rb_{r+1}\ds b_{s}\;|\;1\le r\le s\}$. Note that $P_v$ is a downset 
in $(U,\prec)$. The next theorem was proved by Morse and Hedlund \cite{mors_hedl}, 
see also Allouche and Shallit \cite[Theorem 10.2.6]{allo_shal}.

\begin{theorem}[Morse and Hedlund]\label{subwords}
Let $P$ be the set of all finite subwords of an infinite word $v$ over a finite 
alphabet $A$. Then $f_P(n)=\#\{u\in P\;|\;|u|=n\}$ is either 
larger than $n$ for every $n$ in $\N$ or is eventually constant. In the latter 
case the word $v$ is eventually periodic.
\end{theorem}

\noindent
The case when $P$ is a general downset in $(U,\prec)$, not necessarily coming 
from an infinite word (cf. subsection~\ref{profily}), is discussed below in 
Theorem~\ref{BBwords}.

Examples 1 and 2 are exact results and example 3 combines a 
tight form of an asymptotic inequality with an exact result. Examples 1 and 2 involve only 
countably many counting functions $f_P(n)$ and, as follows from the proofs,
even only countably many downsets $P$. In example 3 we have uncountably many 
distinct counting functions. To see this, take $A=\{0,1\}$ and consider infinite words 
$v$ of the form $v=10^{n_1}10^{n_2}10^{n_3}1\ds$ where $1\le n_1<n_2<n_2<\ds$ is a sequence of integers and $0^m=00\ds 0$ with $m$ zeros. It follows that for distinct 
words $v$ the counting functions $f_{P_v}$ are distinct; Proposition~\ref{uncmany} presents similar arguments in more general settings.

\subsection{Content of the overview}

The previous three examples illuminated to some extent general enumerative results 
we are interested in but they are not fully representative because 
we shall cover a larger area than the growth of downsets. We do not attempt to set 
forth any more formalized definition of a general enumerative result than 
the initial scheme but in subsections~\ref{profily} and \ref{logika} we will discuss some general approaches of finite model theory based on relational structures. 
Not every result or problem mentioned
here fits naturally the scheme; Proposition~\ref{uncmany} and Theorem~\ref{ALthm} are rather 
results to the effect that $\{f_P\;|\;P\in{\cal C}\}$ is too big to be 
contained in a small class ${\cal F}$. This collection 
of general enumerative results is naturally limited by the author's research area and 
his taste but we do hope that it will be of interest to others and that it will inspire a
quest for further generalizations, strengthenings, refinements, common links,
unifications etc.     

For the lack of space, time and expertise we do not mention results on growth in 
algebraic structures, especially the continent of growth in groups; we refer the reader
for information to de la Harpe \cite{harp} (and also to Cameron \cite{came_90}). 
Also, this is not a survey on the class of problems $\#$P in computational complexity theory (see Papadimitriou 
\cite[Chapter 18]{papa}). There are other 
areas of general enumeration not mentioned properly here, for example 0-1 laws (see 
Burris \cite{burr} and Spencer \cite{spen}).
Another reason for omissions of nice general results which should be mentioned here is
simply the author's ignorance---all suggestions, comments and information will be greatly appreciated.
   
In the next subsection we review some notions and definitions from combinatorial enumeration,
in particular we recall the notion of Wilfian formula (polynomial-time counting algorithm).
In Section 2 we review results on growth of downsets in posets of combinatorial 
structures. Subsection 2.1 is devoted to pattern avoiding permutations, Subsections 2.2 and 
2.3 to graphs and related structures, and Subsection 2.4 to relational structures. 
Most of the results in Subsections 2.2 and 2.3 were found by Balogh and Bollob\'as 
and their coauthors (\cite{balo_boll,balo_boll_morr_hype,balo_boll_morr_orgr,balo_boll_morr_tour,
balo_boll_morr_pose,balo_boll_saks_sos,balo_boll_simo,balo_boll_wein_spee,
balo_boll_wein_penu,balo_boll_wein_meas,balo_boll_wein_bell}). We recommend the comprehensive
survey of Bollob\'as \cite{boll_bcc07} on this topic. In Section 3 we advertise five 
topics in general enumeration together with some related results. 
1. The Ehrhart--Macdonald theorem on numbers of lattice points in lattice 
polytopes.  
2. Growth of context-free languages. 3. The theorem of Gessel on  
numbers of labeled regular graphs. 4. The theorem of Bell, Burris and Yeats on 
frequent occurrence of the asymptotics $cn^{-3/2}r^n$. 5. The Specker--Blatter 
theorem on numbers of MSOL-definable structures.

\subsection{Notation and some specific counting functions}

As above, we write $\N$ for the set $\{1,2,3,\dots\}$, $\No$ for $\{0,1,2,\dots\}$ and $[n]$ for 
$\{1,2,\ds,n\}$. We use 
$\#X$ and $|X|$ to denote the cardinality of a set. By the phrase ``for every $n$'' we 
mean ``for every $n$ in $\N$'' and by ``for large 
$n$'' we mean ``for every $n$ in $\N$ with possibly finitely many exceptions''. 
Asymptotic relations are always based on $n\to\infty$. The {\em growth constant} 
$c=c(P)$ of a problem $P$ is $c=\limsup f_P(n)^{1/n}$; the reciprocal 
$\frac{1}{c}$ is then the radius of convergence of the power series $\sum_{n\ge 0}f_P(n)x^n$. 

We review several counting sequences appearing in the mentioned
results. {\em Fibonacci numbers} $(F_n)=(1,2,3,5,8,13,\ds)$ are given by the recurrence 
$F_0=F_1=1$ and $F_n=F_{n-1}+F_{n-2}$ for $n\ge 2$. They are a particular case 
$F_n=F_{n,2}$ of the {\em generalized Fibonacci numbers} $F_{n,k}$, given by
the recurrence $F_{n,k}=0$ for $n<0$, $F_{0,k}=1$ and 
$F_{n,k}=F_{n-1,k}+F_{n-2,k}+\ds+F_{n-k,k}$ for $n>0$. Using the notation 
$[x^n] G(x)$ for the coefficient of $x^n$ in the power series expansion of the expression 
$G(x)$, we have 
$$
F_{n,k}=[x^n]\frac{1}{1-x-x^2-\ds-x^k}. 
$$
Standard methods provide asymptotic relations $F_{n,2}\sim c_2(1.618\ds)^n$, 
$F_{n,3}\sim c_3(1.839\ds)^n$, $F_{n,4}\sim c_4(1.927\ds)^n$ and generally 
$F_{n,k}\sim c_k\alpha_k^n$ for constants $c_k>0$ and $1<\alpha_k<2$; $\frac{1}{\alpha_k}$ 
is the least positive root of the denominator $1-x-x^2-\ds-x^k$ and 
$\alpha_2,\alpha_3,\ds$ monotonicly increase to $2$.
The unlabeled exponential growth of tournaments (Theorem~\ref{BBMtourn}) is governed 
by the {\em quasi-Fibonacci numbers} $F_n^*$ defined by the recurrence $F_0^*=F_1^*=F_2^*=1$ and $F_n^*=F_{n-1}^*+F_{n-3}^*$ for $n\ge 3$; so 
$$
F_n^*=[x^n]\frac{1}{1-x-x^3}
$$
and $F_n^*\sim c(1.466\ds)^n$. 

We introduced {\em Stirling numbers} $S(n,k)$ in example 1. {\em Bell numbers} 
$B_n=\sum_{k=1}^n S(n,k)$ count all partitions of an $n$-elements set and follow
the recurrence $B_0=1$ and $B_n=\sum_{k=0}^{n-1}{n-1\choose k}B_{k}$ for $n\ge 1$. Equivalently,
$$
B_n=[x^n]\sum_{k=0}^{\infty}\frac{x^k}{(1-x)(1-2x)\ds(1-kx)}.
$$
The asymptotic form of the Bell numbers is 
$$
B_n=n^{n(1-\log\log n/\log n+O(1/\log n))}.
$$ 
The numbers $p_n$ of {\em integer partitions} of $n$ count the ways to express $n$ 
as a sum of possibly repeated summands from $\N$, with the order of summands being irrelevant. Equivalently, 
$$
p_n=[x^n]\prod_{k=1}^{\infty}\frac{1}{1-x^k}.
$$
The asymptotic form of $p_n$ is $p_n\sim cn^{-1}\exp(d\sqrt{n})$ for some constants $c,d>0$. 
See Andrews \cite{andr} for more information on these asymptotics and for recurrences satisfied by $p_n$. 

A sequence $f:\;\N\to\C$ is a {\em quasipolynomial} if for every $n$ we have 
$f(n)=a_k(n)n^k+\ds+a_1(n)n+a_0(n)$ where $a_i:\;\N\to\C$ are periodic functions. Equivalently,
$$
f(n)=[x^n]\frac{p(x)}{(1-x)(1-x^2)\ds(1-x^l)}
$$
for some $l$ in $\N$ and a polynomial $p\in\C[x]$. We say that the sequence $f$ is 
{\em holonomic} (other terms are P-recursive and D-finite) if it satisfies for every $n$ (equivalently, for large $n$) a recurrence 
$$
p_k(n)f(n+k)+p_{k-1}(n)f(n+k-1)+\ds+p_0(n)f(n)=0
$$
with polynomial coefficients $p_i\in\C[x]$, not all zero. Equivalently, the power series 
$\sum_{n\ge 0}f(n)x^n$ satisfies a linear differential equation with polynomial coefficients. 
Holonomic sequences generalize sequences satisfying linear recurrences with constant 
coefficients.
The sequences $S(n,k)$, $F_{n,k}$, and $F_n^*$ for each fixed $k$ satisfy a linear recurrence with constant 
coefficients and are holonomic. 
The sequences of Catalan numbers $\frac{1}{n+1}{2n\choose n}$ and of factorial numbers $n!$ are holonomic too. The sequences $B_n$ and $p_n$ are not 
holonomic (\cite{stan2}). It is not hard to show that if $(a_n)$ is holonomic and 
every $a_n$ is in $\Q$, then the polynomials $p_i(x)$ in the recurrence can be taken with 
integer coefficients. In particular, there are only countably many holonomic rational 
sequences.

Recall that a power series $F=\sum_{n\ge 0}a_nx^n$ with $a_n$ in $\C$ is {\em algebraic} 
if there exists a nonzero polynomial $Q(x,y)$ in $\C[x,y]$ such that $Q(x,F(x))=0$. $F$ is 
{\em rational} if $Q$ has degree $1$ in $y$, that is, $F(x)=R(x)/S(x)$ for two polynomials 
in $\C[x]$ where $S(0)\ne 0$. It is well known (Comtet \cite{comt}, Stanley 
\cite{stan2}) that algebraic power series have holonomic coefficients and that the coefficients of rational power series 
satisfy (for large $n$) linear recurrence with constant coefficients. 

{\bf Wilfian formulas.} A counting function $f_P(n)$ has a 
{\em Wilfian formula} (Wilf \cite{wilf_what}) if there exists an algorithm that 
calculates $f_P(n)$ for every input $n$ effectively, that is to say, in polynomial time. 
More precisely, we require (extending the definition in \cite{wilf_what}) 
that the algorithm calculates $f_P(n)$ in the number of steps polynomial in the quantity
$$
t=\max(\log n,\log f_P(n))
$$
---this is (roughly) the minimum time needed for reading the input and writing down 
the answer. In the most common 
situations when $\exp(n^c)<f_P(n)<\exp(n^d)$ for large $n$ and some constants $d>c>0$, this amounts to requiring a number of steps polynomial 
in $n$. But if $f_P(n)$ is small (say $\log n$) or big (say doubly exponential in 
$n$), then one has to work with $t$ in place of $n$. The class of counting functions with 
Wilfian formulas includes holonomic sequences but is much more comprehensive than that.  

\section{Growth of downsets of combinatorial structures}

We survey results in the already 
introduced setting of downsets in posets of combinatorial structures $(U,\prec)$. The function 
$f_P(n)$ counts structures of size $n$ in the downset $P$ and $P$ can also be defined in terms of forbidden substructures as 
$P=\av(F)$. Besides the containment relation $\prec$ we 
employ also isomorphism equivalence relation $\sim$ on $U$ and will count {\em unlabeled}  (i.e., nonisomorphic) structures in $P$. We denote the corresponding counting 
function $g_P(n)$, so
$$
g_P(n)=\#(\{A\in P\;|\;s(A)=n\}/\!\sim)
$$
is the number of isomorphism classes of structures with size $n$ in $P$.

Restrictions on $f_P(n)$ and $g_P(n)$
defining the classes of functions ${\cal F}$ often have the form of {\em jumps in growth}.
A jump is a region of growth prohibited for counting functions---every 
counting function resides either below it or above it. There are many kinds of jumps but 
the most spectacular is perhaps the {\em polynomial--exponential jump} from polynomial to exponential growth, which prohibits counting functions satisfying $n^k<f_P(n)<c^n$ for 
large $n$ for any constants $k>0$ and $c>1$. For groups, Grigorchuk constructed 
a finitely generated group having such intermediate growth (Grigorchuk \cite{grig}, 
Grigorchuk and Pak \cite{grig_pak}, \cite{harp}), which excludes the polynomial--exponential 
jump for general finitely generated groups, but a conjecture says that this jump occurs for 
every finitely presented group. We have seen this jump in Theorems~\ref{setpartstir} and 
\ref{herprostir} (from polynomial growth to growth at least $2^n$) and will meet new examples in 
Theorems \ref{KKthm}, \ref{BBMorgr}, \ref{Kthm}, \ref{BBMtourn}, and \ref{Tthm}.        

If $(U,\prec)$ has an infinite 
antichain $A$, then under natural conditions we 
get uncountably many functions $f_P(n)$. This was observed several times 
in the context of permutation containment and for completeness we give the
argument here again. These natural conditions, which will
always be satisfied in our examples, are {\em finiteness}, for every $n$ there 
are finitely many structures with size $n$ in $U$, and {\em monotonicity}, 
$s(G)\ge s(H)\;\&\;G\prec H$ implies $G=H$ for every $G, H$ in $U$. (Recall that 
$G\prec G$ for every $G$.)

\begin{proposition}\label{uncmany}
If $(U,\prec)$ and the size function $s(\cdot)$ satisfy the monotonicity and 
finiteness conditions and $(U,\prec)$ has an 
infinite antichain $A$, then the set of counting functions $f_P(n)$ is uncountable.
\end{proposition}
\duk
By the assumption on $U$ we can assume that the members of $A$ have distinct sizes. 
We show that all the counting functions $f_{\av(F)}$ for $F \subset A$ are distinct and so this set of functions is uncountable. 
We write simply $f_F$ instead of $f_{\av(F)}$. If $X, Y$ are 
two distinct subsets of $A$, we express them as $X=T\cup\{G\}\cup U$ and 
$Y=T\cup\{H\}\cup V$ so that, without loss of generality, $m=s(G)<s(H)$, and $G_1\in T, G_2\in U$ implies 
$s(G_1)<s(G)<s(G_2)$ and similarly for $Y$ (the sets $T,U,V$ may be empty). Then, by the assumption on $\prec$ and $s(\cdot)$,
$$
f_X(m)=f_{T\cup\{G\}}(m)=f_{T}(m)-1=f_{T\cup\{H\}\cup V}(m)-1=f_{Y}(m)-1
$$
and $f_X\ne f_Y$. 
\kduk

\noindent
An infinite antichain thus gives not only uncountably many downsets but
in fact uncountably many counting functions. Then, in particular, almost all counting functions 
are not computable because we have only countably many algorithms. Recently, Albert 
and Linton \cite{albe_lint} significantly refined this argument by showing how certain infinite 
antichains of permutations produce even uncountably many growth constants, 
see Theorem~\ref{ALthm}. 

On the other hand, if every antichain is finite then there are only countably many 
functions $f_P(n)$. Posets with no infinite antichain are called {\em well quasiorderings} or 
shortly {\em wqo}. (The second part of the wqo property, nonexistence of infinite 
strictly descending chains, is satisfied automatically by the monotonicity condition.) 
But even if $(U,\prec)$ has infinite antichains, there still may be only countably many
downsets $P$ with slow growth of $f_P(n)$. For example, this is the case in 
Theorems~\ref{setpartstir} and \ref{herprostir}. It is then of interest to determine
for which growth uncountably many downsets appear (cf. Theorem~\ref{Vthm}). The posets $(U,\prec)$ considered here usually have infinite antichains, with two notable wqo exceptions consisting of the minor
ordering on graphs and the subsequence ordering on words over a finite alphabet.

\subsection{Permutations}

$U$ is the universe of permutations represented by finite sequences 
$b_1b_2\dots b_n$ such that $\{b_1,b_2,\dots,b_n\}=[n]$. The size of a permutation 
$\pi=a_1a_2\ds a_m$ is its length $|\pi|=m$. The containment relation 
on $U$ is defined by $\pi=a_1a_2\ds a_m\prec\rho=b_1b_2\dots b_n$ if and only if
for some increasing injection $f:\;[m]\to[n]$ one has $a_r<a_s\iff b_{f(r)}<b_{f(s)}$
for every $r,s$ in $[m]$. Problems $P$ are downsets in 
$(U,\prec)$ and their counting functions are
$f_P(n)=\#\{\pi\in P\;|\;|\pi|=n\}$. The poset of permutations
$(U,\prec)$ has infinite antichains (see Spielman and B\'ona \cite{spie_bona}).
For further information and background on the enumeration of downsets of 
permutations see B\'ona \cite{bona_kni}. 

Recall that $c(P)=\limsup f_P(n)^{1/n}$. We define
$$
E=\{c(P)\in[0,+\infty]\;|\;\mbox{$P$ is a downset of permutations}\}
$$
to be the set of growth constants of downsets of permutations. $E$ contains elements
$0,1$ and $+\infty$ because of the downsets $\emptyset$, $\{(1,2,\ds,n)\;|\;n\in\N\}$
and $U$ (all permutations), respectively. How much does $f_P(n)$ drop from $f_U(n)=n!$ if 
$P\ne U$? The {\em Stanley--Wilf conjecture} (B\'ona \cite{bona_97,bona_kni}) asserted that 
it drops to exponential growth. The conjecture was proved in 2004 by Marcus 
and Tardos \cite{marc_tard}.  
\begin{theorem}[Marcus and Tardos]\label{MTthm}
If $P$ is a downset of permutations that is not equal to the set of all permutations,
then, for some constant $c$, $f_P(n)<c^n$ for every $n$.
\end{theorem}

\noindent 
%konstanty, Cibulka.
Thus, with the sole exception of $U$, every $P$ has a finite 
growth constant. Arratia \cite{arra} showed that if $F=\{\pi\}$ then  
$c(P)=c(\av(F))$ is attained as a limit $\lim f_P(n)^{1/n}$. It would be nice to 
extend this result.

\begin{problem}
Does $\lim f_P(n)^{1/n}$ always exist when $F$ in $P=\av(F)$ has more than one forbidden 
permutation?
\end{problem}

\noindent
For infinite $F$ there conceivably might be oscillations between two different 
exponential growths (similar oscillations occur for hereditary graph properties and for downsets 
of words). It would be surprising if oscillations occurred for finite $F$.

Kaiser and Klazar \cite{kais_klaz} determined growths of downsets of permutations 
in the range up to $2^{n-1}$.  

\begin{theorem}[Kaiser and Klazar]\label{KKthm}
If $P$ is a downset of permutations, then exactly one of the four cases occurs. 
\begin{enumerate}
\item For large $n$, $f_P(n)$ is constant.
\item There are integers $a_0,\ds,a_k$, $k\ge 1$ and $a_k>0$, such that  
$f_P(n)=a_0{n\choose 0}+\ds+a_k{n\choose k}$ for large $n$. 
Moreover, $f_P(n)\ge n$ for every $n$. 
\item There are constants $c,k$ in $\N$, $k\ge 2$, such that 
$F_{n,k}\le f_P(n)\le n^cF_{n,k}$ for every $n$, where $F_{n,k}$ are the generalized 
Fibonacci numbers. 
\item One has $f_P(n)\ge 2^{n-1}$ for every $n$.
\end{enumerate}
The lower bounds in cases 2, 3, and 4 are best possible.
\end{theorem}

\noindent
This implies that
$$
E\cap[0,2]=\{0,1,2,\alpha_2,\alpha_3,\alpha_4,\ds\},
$$
$\alpha_k$ being the growth constants of $F_{n,k}$, 
and that $\lim f_P(n)^{1/n}$ exists and equals to $0,1$ or to some $\alpha_k$ whenever 
$f_P(n)<2^{n-1}$ for one $n$. Note that $2$ is the single accumulation point of $E\cap[0,2]$.
We shall see that Theorem~\ref{KKthm} is subsumed in Theorem~\ref{BBMorgr} on ordered graphs.
Case 2 and case 3 with $k=2$ give the {\em polynomial--Fibonacci jump}: 
If  $P$ is a downset of permutations, then either $f_P(n)$ grows at most polynomially 
(and in fact equals  to a polynomial for large $n$) or at least 
Fibonaccially. Huczynska and Vatter \cite{hucz_vatt} gave a simpler proof 
for this jump. Theorem~\ref{Kthm} extends it to edge-colored cliques. 
Theorem~\ref{KKthm} combines an exact result in case 1 and 2 with an asymptotic result in 
case 3. It would be nice to have in case 3 an exact result too and to determine  
precise forms of the corresponding functions $f_P(n)$ (it is known that in cases 1--3 
the generating function $\sum_{n\ge 0}f_P(n)x^n$ is rational, see the remarks at the end of 
this subsection). Klazar \cite{klaz_thcs} proved that 
cases 1--3 comprise only countably many downsets, more precisely: if
$f_P(n)<2^{n-1}$ for one $n$, then $P=\av(F)$ has finite base $F$.
In the other direction he showed (\cite{klaz_thcs}) that there are uncountably many downsets 
$P$ with $f_P(n)<(2.336\ds)^n$ for large $n$. Recently, Vatter \cite{vatt_manu} 
determined the uncountability threshold precisely and extended the description of $E$ above 
$2$. 

\begin{theorem}[Vatter]\label{Vthm}
Let $\kappa=2.205\ds$ be the real root of $x^3-2x^2-1$. There are uncountably 
many downsets of permutations $P$ with $c(P)\le\kappa$ but only countably many of them have 
$c(P)<\kappa$ and for each of these $\lim f_P(n)^{1/n}$ exists. Moreover, the countable intersection
$$
E\cap(2,\kappa)
$$
consists exactly of the largest positive roots of the polynomials in the four families
($k,l$ range over $\N$)
\begin{enumerate}
\item $3-x-x^{k+1}-2x^{k+3}+x^{k+4}$,  
\item $1+2x-x^2-x^{k+2}-2x^{k+4}+x^{k+5}$, 
\item $1+x^k-x^{k+l}-2x^{k+l+2}+x^{k+l+3}$, and 
\item $1-x^k-2x^{k+2}+x^{k+3}$. 
\end{enumerate}
\end{theorem}

\noindent
The set $E\cap(2,\kappa)$ has no accumulation point from above but it has infinitely 
many accumulation points from below: $\kappa$ is the smallest element of $E$ which is an accumulation point of accumulation points. The smallest element of $E\cap(2,\kappa)$ is 
$2.065\ds$ ($k=l=1$ in the family 3).

In \cite{balo_boll_morr_orgr} it was conjectured that all elements of $E$ 
(even in the more general situation of ordered graphs) are algebraic numbers and that $E$ 
has no accumulation point from above. These conjectures were refuted by Albert and Linton
\cite{albe_lint}. Recall that a subset of $\R$ is {\em perfect} if it is closed and has no 
isolated point. Due to the completeness of $\R$ such a set is inevitably uncountable. 

\begin{theorem}[Albert and Linton]\label{ALthm}
The set $E$ of growth constants of downsets of permutations contains a perfect
subset and therefore is uncountable. Also, $E$ contains accumulation points from above. 
\end{theorem}

\noindent
The perfect subset constructed by Albert and Linton has smallest element $2.476\ds$ and  
they conjecture that $E$ contains some real interval $(\lambda,+\infty)$. However, a
typical downset produced by their construction has infinite base. It seems that the 
refuted conjectures should have been phrased for finitely based downsets.

\begin{problem}  
Let $E^*$ be the countable subset of $E$ consisting of the growth constants
of finitely based downsets of permutations. Show that every $\alpha$ in $E^*$ 
is an algebraic number and that for every $\alpha$ in $E^*$ there is a $\delta>0$ such that 
$(\alpha,\alpha+\delta)\cap E^*=\emptyset$.  
\end{problem}

\noindent
We know that $E^*\cap[0,2]=E\cap[0,2]$ and probably (as conjectured in \cite{vatt_manu}) 
even $E^*\cap[0,\kappa)=E\cap[0,\kappa)$.

We turn to the questions of exact counting. In view of Proposition~\ref{uncmany} and 
Theorem~\ref{ALthm}, we restrict to downsets of permutations with finite bases.  
The next problem goes back to Gessel \cite[the final section]{gess}. 

\begin{problem}
Is it true that for every finite set of permutations $F$ the counting function 
$f_{\av(F)}(n)$ is holonomic? 
\end{problem}

\noindent
All explicit $f_P(n)$ found so far are holonomic. Zeilberger conjectures (\cite{elde_vatt}) 
that $P=\av(1324)$ has nonholonomic counting function (see Marinov and Radoi\'ci\v c 
\cite{mari_rado} and Albert et al. \cite{albe_al} for the approaches to counting $\av(1324)$). 
We remarked earlier that almost all infinitely based $P$ have nonholonomic $f_P(n)$.

More generally, one may pose (Vatter \cite{vatt_schemes}) the following 
question. 

\begin{problem}
Is it true that for every finite set of permutations $F$ the counting function 
$f_{\av(F)}(n)$ has a Wilfian formula, that is, can be evaluated by an algorithm in number 
of steps polynomial in $n$? 
\end{problem}

\noindent
Wilfian formulas were shown to exist for several classes of finitely based downsets of permutations. We refer the reader to Vatter \cite{vatt_schemes} for further information and 
mention here only one such result due to Albert 
and Atkinson \cite{albe_atki}. Recall that $\pi=a_1a_2\ds a_n$ is a {\em simple permutation}
if $\{a_i,a_{i+1},\ds,a_j\}$ is not an interval in $[n]$ for every $1\le i\le j\le n$, 
$0<j-i<n-1$.

\begin{theorem}[Albert and Atkinson]\label{AAthm}
If $P$ is a downset of permutations containing only finitely many simple permutations, 
then $P$ is finitely based and the generating function $\sum_{n\ge 0}f_P(n)x^n$ is 
algebraic and thus $f_P(n)$ has a Wilfian formula.
\end{theorem}

\noindent
Brignall, Ru\v skuc and Vatter \cite{brig_rusk_vatt} show that it is decidable whether a 
downset given by its finite basis contains finitely many simple permutations 
and Brignall, Huczynska and Vatter \cite{brig_hucz_vatt} extend Theorem~\ref{AAthm} by showing
that many subsets of downsets with finitely many simple permutations have algebraic
generating functions as well.

We conclude this subsection by looking back at Theorems \ref{KKthm} and
\ref{Vthm} from the standpoint of effectivity. Let a downset of permutations $P=\av(F)$
be given by its finite base $F$. Then it is decidable whether $c(P)<2$ and (as noted in 
\cite{vatt_manu}) for $c(P)<2$ the results of Albert, Linton and Ru\v skuc \cite{albe_lint_rusk} provide effectively a Wilfian formula for $f_P(n)$, in fact, the generating function
is rational. Also, it is decidable whether $f_P(n)$ is a polynomial for large $n$ 
(\cite{hucz_vatt}, Albert, Atkinson and Brignall \cite{albe_atki_brig}). By \cite{vatt_manu}, 
it is decidable whether $c(P)<\kappa$ and Vatter conjectures that even for $c(P)<\kappa$ the 
generating function of $P$ is rational. 

%Bonovy vysledky o Av(pi)

\subsection{Unordered graphs}

$U$ is the universe of finite simple graphs with normalized vertex
sets $[n]$ and $\prec$ is the induced subgraph relation. 
Problems $P$ are hereditary graph properties, that is, downsets in $(U,\prec)$, and 
$f_P(n)$ counts the graphs in $P$ with $n$ vertices. 
A more restricted family is {\em monotone properties}, which are 
hereditary properties that are closed under taking any subgraph. An even more restricted
family consists of {\em minor-closed classes}, which are monotone 
properties that are closed under contracting edges. By Proposition~\ref{uncmany} there are uncountably many 
counting functions of monotone properties (and hence of hereditary properties as well) 
because, for example, the set of all cycles is an infinite antichain to subgraph
ordering. On the other hand, by the monumental theorem of Robertson and Seymour 
\cite{minors} there are no infinite antichains in the minor ordering and so there are only countably many minor-closed classes. The following remarkable theorem describes 
growths of hereditary properties.

\begin{theorem}[Balogh, Bollob\'as, Weinreich, Alekseev, Thomason]\label{BBWherpro}
If $P$ is a proper hereditary graph property then exactly one of the four cases occurs.
\begin{enumerate}
\item There exist rational polynomials $p_1(x),p_2(x),\dots,p_k(x)$ such that 
$f_P(n)=p_1(n)+p_2(n)2^n+\dots+p_k(n)k^n$ for large $n$. 
\item There is a constant $k$ in $\N$, $k\ge 2$, such that $f_P(n)=n^{(1-1/k)n+o(n)}$ 
for every $n$. 
\item One has $n^{n+o(n)}<f_P(n)<2^{o(n^2)}$ for every $n$. 
\item There is a constant $k$ in $\N$, $k\ge 2$, such that 
$f_P(n)=2^{(1/2-1/2k)n^2+o(n^2)}$ for every $n$.  
\end{enumerate}   
\end{theorem} 

\noindent
We mentioned case 1 as Theorem~\ref{herprostir}. The first three cases were proved
by Balogh, Bollob\'as and Weinreich in \cite{balo_boll_wein_spee}. The fourth case is 
due to Alekseev \cite{alek} and independently Bollob\'as and Thomason 
\cite{boll_thoma_proj}. 

Now we will discuss further strengthenings and refinements of 
Theorem~\ref{BBWherpro}. 
Scheinerman and Zito in a pioneering work \cite{sche_zito} obtained its weaker version. They showed that for a hereditary graph property $P$ 
either (i) $f_P(n)$ is constantly $0$, $1$ or $2$ for large $n$ or (ii) 
$an^k<f_P(n)<bn^k$ for every $n$ and some constants $k$ in $\N$ and $0<a<b$ 
or (iii) $n^{-c}k^n\le f_P(n)\le n^ck^n$ for every $n$ and constants $c,k$ in $\N$, 
$k\ge 2$, or (iv) $n^{cn}\le f_P(n)\le n^{dn}$ for 
every $n$ and some constants $0<c<d$ or (v) $f_P(n)>n^{cn}$ for large $n$ for 
every constant $c>0$. 

In cases 1, 2, and 4 growths of $f_P(n)$ settle to specific asymptotic values 
and these can be characterized by certain minimal hereditary properties, as shown 
in \cite{balo_boll_wein_spee}. Case 3, the penultimate rate of growth 
(\cite{balo_boll_wein_penu}), is very different. Balogh, Bollob\'as and Weinreich 
proved in \cite{balo_boll_wein_penu} that for every $c>1$ and $\varepsilon>1/c$ there 
is a monotone property $P$ such that 
$$
f_P(n)\in[n^{cn+o(n)},\; 2^{(1+o(1))n^{2-\varepsilon}}]
$$
for every $n$ and $f_P(n)$ attains either extremity of the interval
infinitely often. Thus in case 3 the growth may oscillate (infinitely often) between 
the bottom and top parts of the range. The paper \cite{balo_boll_wein_penu} 
contains further examples of oscillations (we stated here just one simplified version) 
and a conjecture that for finite $F$ the functions $f_{\av(F)}(n)$ do not oscillate. 
As for the
upper boundary of the range, in \cite{balo_boll_wein_penu} it is proven that for every
monotone property $P$, 
$$
f_P(n)=2^{o(n^2)}\Rightarrow f_P(n)<2^{n^{2-1/t+o(1)}}\ \mbox{ for some $t$ in $\N$}.
$$
For hereditary properties this jump is only conjectured. What about the lower boundary?
The paper \cite{balo_boll_wein_bell} is devoted to the proof of the following theorem.

\begin{theorem}[Balogh, Bollob\'as and Weinreich]
If $P$ is a hereditary graph property, then exactly one of the two cases occurs.  
\begin{enumerate}
\item There is a constant $k$ in $\N$ such that $f_P(n)<n^{(1-1/k)n+o(n)}$ for every 
$n$.
\item For large $n$, one has $f_P(n)\ge B_n$ where $B_n$ are Bell numbers. This lower
bound is the best possible.
\end{enumerate}
\end{theorem}

\noindent
By the theorem, the growth of Bell numbers is the lower boundary of the penultimate 
growth in case 3 of Theorem~\ref{BBWherpro}. 

Monotone properties of graphs are hereditary and therefore their counting functions follow
Theorem~\ref{BBWherpro}. Their more restricted nature allows for
simpler proofs and simple characterizations of minimal monotone properties, which 
is done in the paper \cite{balo_boll_wein_meas}. Certain growths of hereditary properties
do not occur for monotone properties, for example if $P$ is monotone and $f_P(n)$ is 
unbounded, then $f_P(n)\ge{n\choose 2}+1$ for every $n$ (\cite{balo_boll_wein_meas}) but,
$P$ consisting of complete graphs with possibly an additional isolated vertex is a hereditary 
property with $f_P(n)=n+1$ for $n\ge 3$. More generally, Balogh, Bollob\'as and Weinreich
show in \cite{balo_boll_wein_meas} that if $P$ is monotone and $f_P(n)$ grows polynomially,
then 
$$
f_P(n)=a_0{n\choose 0}+a_1{n\choose 1}+\ds+a_k{n\choose k}\ \mbox{ for large $n$ }\  
$$
and some integer constants $0\le a_j\le 2^{j(j-1)/2}$. In fact, \cite{balo_boll_wein_meas} 
deals mostly with general results on the extremal functions 
$e_P(n):=\max\{|E|\;|\;G=([n],E)\in P\}$ for monotone properties $P$.

For the top growths in case 4 of Theorem~\ref{BBWherpro}, Alekseev \cite{alek} 
and Bollob\'as and Thomason 
\cite{boll_thoma_here} proved that for $P=\av(F)$ with $f_P(n)=2^{(1/2-1/2k)n^2+o(n^2)}$
the parameter $k$ is equal to the maximum $r$ such that there is an $s$, $0\le s\le r$, 
with the property that no graph in $F$ can have its vertex set partitioned into  
$r$ (possibly empty) blocks inducing $s$ complete 
graphs and $r-s$ empty graphs. For monotone properties this reduces
to $k=\min\{\chi(G)-1\;|\;G\in F\}$ where $\chi$ is the usual chromatic number of graphs. 
Balogh, Bollob\'as and Simonovits \cite{balo_boll_simo} replaced for monotone properties the error term $o(n^2)$ by $O(n^{\gamma})$, $\gamma=\gamma(F)<2$. Ishigami \cite{ishi} recently
extended case 4 to $k$-uniform hypergraphs. 

%Btw, kdy se objevuje nespocetne mnoho her. properties? 

Minor-closed classes of graphs were recently looked at from the point of view 
of counting functions as well. They again follow Theorem~\ref{BBWherpro}, with possible simplifications due to their more restricted nature. One is that there are only 
countably many minor-closed classes.
Another simplification is that, with the trivial 
exception of the class of all graphs, case 4 does not occur as proved by 
Norine et al. \cite{nori_all}. 

\begin{theorem}[Norine, Seymour, Thomas and Wollan]
If $P$ is a proper minor-closed class of graphs then $f_P(n)<c^nn!$ for every $n$ in 
$\N$ for a constant $c>1$.
\end{theorem}

Bernardi, Noy and Welsh \cite{bern_noy_wels} obtained the following theorem; 
we shorten its statement by omitting characterizations of classes $P$ with 
the given growth rates.  

\begin{theorem}[Bernardi, Noy and Welsh]
If $P$ is a proper minor-closed class of graphs then exactly one of the six cases 
occurs. 
\begin{enumerate}
\item The counting function $f_P(n)$ is constantly $0$ or $1$ for 
large $n$.     
\item For large $n$, $f_P(n)=p(n)$ for a rational polynomial $p$ of degree
at least $2$. 
\item For every $n$, $2^{n-1}\le f_P(n)<c^n$ for a constant 
$c>2$.
\item There exist constants $k$ in $\N$, $k\ge 2$, and $0<a<b$ such that 
$a^n n^{(1-1/k)n}<f_P(n)<b^n n^{(1-1/k)n}$ 
for every $n$.  
\item For every $n$, $B_n\le f_P(n)=o(1)^nn!$ where $B_n$ is the $n^{\mbox{\scriptsize{th}}}$ Bell number. 
\item For every $n$, $n!\le f_P(n)<c^nn!$ for a constant $c>1$.
\end{enumerate} 
The lower bounds in cases 3, 5 and 6 are best possible. 
\end{theorem}

\noindent
In fact, in case 3 the formulas of Theorem~\ref{herprostir} apply. Using the
strongly restricted nature of minor-closed classes, one could perhaps obtain in case 
3 an even more specific exact result. Paper \cite{bern_noy_wels} gives
further results on the growth constants $\lim\;(f_P(n)/n!)^{1/n}$ in case 6 and states 
several open problems, of which we mention the following analogue of 
Theorem~\ref{MTthm} for unlabeled graphs. Similar conjecture was stated also in 
McDiarmid, Steger and Welsh \cite{mcdi_steg_wels}.

\begin{problem}
Does every proper minor-closed class of graphs contain at most $c^n$ nonisomorphic
graphs on $n$ vertices, for a constant $c>1$?
\end{problem}

This brings us to the unlabeled count of hereditary properties. The following theorem 
was obtained by Balogh et al. \cite{balo_boll_saks_sos}.

\begin{theorem}[Balogh, Bollob\'as, Saks and S\'os]\label{unlabcount}
If $P$ is a hereditary graph property and $g_P(n)$ counts nonisomorphic graphs in $P$ 
by the number of vertices, then exactly one of the three cases occurs.
\begin{enumerate}
\item For large $n$, $g_P(n)$ is constantly $0,1$ or $2$.
\item For every $n$, $g_P(n)=cn^k+O(n^{k-1})$ for some constants $k$ in $\N$ and $c$ in 
$\Q$, $c>0$. 
\item For large $n$, $g_P(n)\ge p_n$ where $p_n$ is the number of integer partitions of 
$n$. This lower bound is best possible.
\end{enumerate}
\end{theorem}

\noindent
(We have shortened the statement by omitting the characterizations of $P$ with 
given growth rates.) The authors of 
\cite{balo_boll_saks_sos} remark that with more effort case 2 
can be strengthened, for large $n$, to an exact result with the error term 
$O(n^{k-1})$ replaced by a quasipolynomial $p(n)$ of degree at most $k-1$.
It turns out that a weaker form of the jump from case 2 to case 3  
was proved already by Macpherson \cite{macp1, macp2}: If 
$G=(\N,E)$ is an infinite graph and $g_G(n)$ is the number of its unlabeled $n$-vertex
induced subgraphs then either $g_G(n)\le n^c$ for every $n$ and a constant $c>0$ 
or $g_G(n)>\exp(n^{1/2-\ep})$ for large $n$ for every constant $\ep>0$.
Pouzet \cite{pouz_81} showed that in the former case $c_1n^d<g_G(n)<c_2n^d$ for every
$n$ and some constants $0<c_1<c_2$ and $d$ in $\N$. 

\subsection{Ordered graphs and hypergraphs, edge-colored cliques, words, posets, tournaments, and tuples}
\label{subsordered}

%Obecne bla, ze nektere struktury maji usporadane nosne mnoziny.

{\bf Ordered graphs. } As previously, $U$ is the universe of finite simple graphs with
vertex sets $[n]$ but $\prec$ is now the {\em ordered} induced subgraph relation, 
which means that $G_1=([m],E_1)\prec G_2=([n],E_2)$ if and only if there is an increasing injection 
$f:\;[m]\to[n]$ such that $\{u,v\}\in E_1\iff\{f(u),f(v)\}\in E_2$. Problems are downsets
$P$ in $(U,\prec)$, are called {\em hereditary properties of ordered graphs}, 
and $f_P(n)$ is the number of graphs in $P$ with vertex set $[n]$. The next theorem, 
proved by Balogh, Bollob\'as and Morris \cite{balo_boll_morr_orgr}, vastly generalizes 
Theorem~\ref{KKthm}.

\begin{theorem}[Balogh, Bollob\'as and Morris]\label{BBMorgr}
If $P$ is a hereditary property of ordered graphs, then exactly one of the four cases occurs. 
\begin{enumerate}
\item For large $n$, $f_P(n)$ is constant.
\item There are integers $a_0,\ds,a_k$, $k\ge 1$ and $a_k>0$, such that  
$f_P(n)=a_0{n\choose 0}+\ds+a_k{n\choose k}$ for large $n$. 
Moreover, $f_P(n)\ge n$ for every $n$. 
\item There are constants $c,k$ in $\N$, $k\ge 2$, such that 
$F_{n,k}\le f_P(n)\le n^cF_{n,k}$ for every $n$, where $F_{n,k}$ are the generalized 
Fibonacci numbers. 
\item One has $f_P(n)\ge 2^{n-1}$ for every $n$.
\end{enumerate}
The lower bounds in cases 2, 3, and 4 are best possible.
\end{theorem}

\noindent
This is an extension of Theorem~\ref{KKthm} because the poset of permutations is embedded
in the poset of ordered graphs via representing a permutation $\pi=a_1a_2\ds a_n$ by 
the graph $G_{\pi}=([n],\{\{i,j\}\;|\;i<j\;\&\;a_i<a_j\})$. One can check that 
$\pi\prec\rho\iff G_{\pi}\prec G_{\rho}$ and that the graphs $G_{\pi}$ form a downset
in the poset of ordered graphs, so Theorem~\ref{BBMorgr} implies Theorem~\ref{KKthm}. 
Similarly, the poset of set partitions of example 1 in Introduction is embedded in the poset of ordered 
graphs, via representing partitions by graphs whose components are cliques. Thus the growths 
of downsets of set partitions in the range up to $2^{n-1}$ are described by 
Theorem~\ref{BBMorgr}. As for permutations, it would be nice to have in case 3 an exact result.
Balogh, Bollob\'as and Morris conjecture that $2.031\ds$ (the largest real root 
of $x^5-x^4-x^3-x^2-2x-1$) is the smallest growth constant for ordered graphs 
above $2$ and Vatter (\cite{vatt_manu}) notes that this is not an element of $E$ 
and thus here the growth constants for permutations and for ordered graphs
part ways.

{\bf Edge-colored cliques.} Klazar \cite{klaz_grrat} considered the universe $U$ of 
pairs $(n,\chi)$ where $n$ ranges over $\N$ and $\chi$ is a mapping from the set 
${[n]\choose 2}$ of two-element subsets of $[n]$ to a finite set of colors $C$. The
containment $\prec$ is defined by $(m,\phi)\prec(n,\chi)$ if and only if there is an increasing 
injection $f:\;[m]\to[n]$ such that $\chi(\{f(x),f(y)\})=\phi(\{x,y\})$ for every $x,y$ in $[m]$. For two colors we recover ordered graphs with induced ordered subgraph relation. In 
\cite{klaz_grrat} the following theorem was proved. 

\begin{theorem}[Klazar]\label{Kthm}
If $P$ is a downset of edge-colored cliques, then exactly one of the three cases occurs.  
\begin{enumerate}
\item The function $f_P(n)$ is constant for large $n$. 
\item There is a constant $c$ in $\N$ such that $n\le f_P(n)\le n^c$ for every $n$.
\item One has $f_P(n)\ge F_n$ for every $n$, where $F_n$ are the Fibonacci numbers.   
\end{enumerate}
The lower bounds in cases 2 and 3 are best possible.
\end{theorem}

\noindent
This extends the bounded-linear jump and the polynomial-Fibonacci jump of 
Theorem~\ref{BBMorgr}. It would be 
interesting to have full Theorem~\ref{BBMorgr} in this more general setting. As 
explained in \cite{klaz_grrat}, many posets of structures can be embedded in the 
poset of edge-colored cliques (as we have just seen for permutations) 
and thus Theorem~\ref{Kthm} applies to them. With more effort, case 2 can be 
strengthened to the exact result $f_P(n)=p(n)$ with rational polynomial $p(x)$. 

{\bf Words over finite alphabet.} We revisit example 3 from Introduction.
Recall that $U=A^*$ consists of all finite words over finite alphabet $A$ and $\prec$
is the subword ordering. This ordering has infinite antichains, for example $11, 101,
1001,\ds$ for $A=\{0,1\}$. Balogh and Bollob\'as \cite{balo_boll} investigated general downsets in $(A^*,\prec)$ and proved the following extension of Theorem~\ref{subwords}.

\begin{theorem}[Balogh and Bollob\'as]\label{BBwords}
If $P$ is a downset of finite words over a finite alphabet $A$ in the subword 
ordering, then $f_P(n)$ is either bounded or $f_P(n)\ge n+1$ for every $n$.  
\end{theorem}

\noindent
In contrast with Theorem~\ref{subwords}, for general downsets, a bounded function  
$f_P(n)$ need not be eventually constant. Balogh and Bollob\'as \cite{balo_boll}
showed that for fixed $s$ in $\N$ function $f_P(n)$ may oscillate infinitely often between the maximum and minimum values $s^2$ and $2s-1$, and $s^2+s$ and $2s$. These 
are, however, the wildest bounded oscillations possible since they proved, 
as their main result, that if $f_P(n)=m\le n$ for some $n$ then $f_P(N)\le(m+1)^2/4$ 
for every $N$, $N\ge n+m$. They also gave examples of unbounded oscillations of $f_P(n)$
between $n+g(n)$ and $2^{n/g(n)}$ for any increasing and unbounded function 
$g(n)=o(\log n)$, with the downset $P$ coming from an infinite word over two-letter 
alphabet.

Another natural ordering on $A^*$ is the {\em subsequence 
ordering}, $a_1a_2\ds a_k\prec b_1b_2\ds b_l$ if and only if 
$b_{i_1}=a_1, b_{i_2}=a_2,\ds,b_{i_k}=a_k$ for some indices $1\le i_1<i_2<\ds<i_k\le l$.
Downsets in this ordering remain downsets in the subword ordering and 
thus their counting functions are governed by Theorem~\ref{BBwords}. But they can be also embedded in the poset of edge-colored complete graphs (associate with $a_1a_2\ds a_n$ 
the pair $(n,\chi)$ where $\chi(\{i,j\})=\{a_i,a_j\}$ for $i<j$) and 
Theorem~\ref{Kthm} applies. In particular, if $P\subset A^*$ is a downset in the 
subsequence ordering, then $f_P(n)$ is constant for large $n$ or $f_P(n)\ge n+1$ for 
every $n$ (by Theorems~\ref{Kthm} and \ref{BBwords}). The subsequence ordering on $A^*$ is a
wqo by Higman's theorem (\cite{higm}) and therefore has only countably many downsets.

A variation on the subsequence ordering is the ordering on $A^*$ given by 
$u=a_1a_2\ds a_k\prec v=b_1b_2\ds b_l$ if and only if there is a permutation $\pi$ of the alphabet $A$
such that $a_1a_2\ds a_k\prec \pi(b_1)\pi(b_2)\ds \pi(b_l)$ in the subsequence ordering,
that is, $u$ becomes a subsequence of $v$ after the letters in $v$ are injectively renamed.
This ordering on $A^*$ gives example 1 in Introduction and leads to 
Theorem~\ref{setpartstir}. It is a wqo as well.

{\bf Posets and tournaments.} $U$ is the set of all pairs 
$S=([n],\le_S)$ where $\le_S$ is a non-strict partial ordering on $[n]$. 
We set $R=([m],\le_R)\prec S=([n],\le_S)$ if and only if there is an injection 
$f:\;[m]\to[n]$ such that $x\le_Ry\iff f(x)\le_Sf(y)$ for every $x,y$ in $[m]$. Thus 
$R\prec S$ means that the poset $R$ is an induced subposet of $S$. Downsets in
$(U,\prec)$, {\em hereditary properties of posets}, and their growths were investigated 
by Balogh, Bollob\'as and Morris in \cite{balo_boll_morr_pose}. For the unlabeled 
count they obtained the following result.

\begin{theorem}[Balogh, Bollob\'as and Morris]
If $P$ is a hereditary property of posets and $g_P(n)$ counts nonisomorphic posets in 
$P$ by the number of vertices, then exactly one of the three cases occurs.
\begin{enumerate}
\item Function $g_P(n)$ is bounded.
\item There is a constant $c>0$ such that, for every $n$, 
$\left\lceil\frac{n+1}{2}\right\rceil\le g_P(n)\le\left\lceil\frac{n+1}{2}\right\rceil+c$. 
\item For every $n$, $g_P(n)\ge n$.
\end{enumerate}
The lower bounds in cases 2 and 3 are best possible.
\end{theorem}

\noindent
As for the labeled count $f_P(n)$, using case 1 of Theorem~\ref{BBWherpro} they proved (\cite[Theorem 2]{balo_boll_morr_pose}) that if $P$ is a hereditary property of posets
then either (i) $f_P(n)$ is constantly $1$ for large $n$ or (ii) there are $k$ 
integers $a_1,\ds,a_k$, $a_k\ne 0$, such that 
$f_P(n)=a_0{n\choose 0}+\ds+a_k{n\choose k}$ for large $n$ or (iii) 
$f_P(n)\ge 2^n-1$ for every $n$, $n\ge 6$. Moreover, the lower bound in case (iii) is 
best possible and in case (ii) one has 
$f_P(n)\ge {n\choose 0}+\ds+{n\choose k}$ for every $n$, $n\ge 2k+1$ and this bound is 
also best possible.

A {\em tournament} is a pair $T=([n],T)$ where $T$ is a binary relation on $[n]$ such that 
$xTx$ for no $x$ in $[n]$ and for every two distinct elements $x,y$ in $[n]$ exactly one
of $xTy$ and $yTx$ holds. $U$ consists of all tournaments for $n$ ranging in $\N$ and 
$\prec$ is the induced subtournament relation. Balogh, 
Bollob\'as and Morris considered  in \cite{balo_boll_morr_tour,balo_boll_morr_pose} unlabeled counting functions of hereditary properties 
of tournaments. We merge their results in a single theorem.

\begin{theorem}[Balogh, Bollob\'as and Morris]\label{BBMtourn}
If $P$ is a hereditary property of tournaments and $g_P(n)$ counts nonisomorphic tournaments in $P$ by the number of vertices, then exactly one of the three cases 
occurs.
\begin{enumerate}
\item For large $n$, function $g_P(n)$ is constant. 
\item There are constants $k$ in $\N$ and $0<c<d$ such that $cn^k<g_P(n)<dn^k$ for 
every $n$. Moreover, $g_P(n)\ge n-2$ for every $n$, $n\ge 4$.
\item For every $n$, $n\ne 4$, one has $g_P(n)\ge F_n^*$ where $F_n^*$ are the
quasi-Fibonacci numbers.
\end{enumerate}
The lower bounds in cases 2 and 3 are best possible.
\end{theorem}

\noindent
Case 1 and the second part of case 2 were proved in \cite{balo_boll_morr_pose} and the rest of the theorem in \cite{balo_boll_morr_tour}. A closely related and in one direction stronger theorem was independently obtained by Boudabbous and Pouzet \cite{boud_pouz} 
(see also \cite[Theorem 22]{pouz_surv}): If $g_T(n)$ counts unlabeled $n$-vertex subtournaments of an infinite 
tournament $T$, then either $g_T(n)$ is a quasipolynomial for large $n$ or 
$g_T(n)>c^n$ for large $n$ for a constant $c>1$. 

{\bf Ordered hypergraphs.} $U$ consists of 
all hypergraphs, which are the pairs $H=([n],H)$ with $n$ in $\N$ and $H$ being a set of 
nonempty and non-singleton subsets of $[n]$, called edges. Note that $U$ extends both
the universe of finite simple graphs and the universe of set partitions. 
The containment $\prec$ is ordered but non-induced and is defined by $([m],G)\prec([n],H)$ if and only if there is an increasing injection $f:\;[m]\to[n]$ and an injection $g:\;G\to H$ such that for every edge $E$ in $G$ we have 
$f(E)\subset g(E)$. Equivalently, one can omit some vertices from $[n]$, some edges 
from $H$ and delete some vertices from the remaining edges in $H$ so that the resulting hypergraph is order-isomorphic to $G$. Downsets in $(U,\prec)$ are called 
{\em strongly monotone properties of ordered hypergraphs}. Again, $f_P(n)$ counts hypergraphs in $P$ with the vertex set $[n]$. We define a special downset $\Pi$: we 
associate with every permutation $\pi=a_1a_2\ds a_n$ the (hyper)graph 
$G_{\pi}=([2n],\{\{i,n+a_i\}\;|\;i\in[n]\})$ and let $\Pi$ denote the set of all hypergraphs in $U$ contained in some 
graph $G_{\pi}$; the graphs in $\Pi$ differ from $G_{\pi}$'s only in 
adding in all ways isolated vertices. Note that $\pi\prec\rho$ for two permutations if and only if 
$G_{\pi}\prec G_{\rho}$ for the corresponding (hyper)graphs. The next theorem 
was conjectured by Klazar in \cite{klaz_stir} for set partitions and in \cite{klaz_nch} 
for ordered hypergraphs. 

\begin{theorem}[Balogh, Bollob\'as, Morris, Klazar, Marcus]\label{BBMKMthm}
If $P$ is a strongly monotone property of ordered hypergraphs, then exactly one of the 
two cases occurs.
\begin{enumerate}
\item There is a constant $c>1$ such that $f_P(n)\le c^n$ for every $n$.
\item One has $P\supset\Pi$, which implies that
$$
f_P(n)\ge\sum_{k=0}^{\lfloor n/2\rfloor}{n\choose 2k}k!=n^{n+O(n/\log n)}
$$
for every $n$ and that the lower bound is best possible.
\end{enumerate}
\end{theorem}

\noindent
The theorem was proved by Balogh, Bollob\'as and Morris 
\cite{balo_boll_morr_hype} and independently by Klazar and Marcus \cite{klaz_marc} 
(by means of results from \cite{klaz_DS, klaz_nch, marc_tard}).  It follows that 
Theorem~\ref{BBMKMthm} implies Theorem~\ref{MTthm}. Theorem~\ref{BBMKMthm} was 
motivated by efforts to extend the Stanley--Wilf conjecture, now the 
Marcus--Tardos theorem, from permutations to more general structures. 
A further extension would be to have it 
for the wider class of hereditary properties of ordered hypergraphs. These correspond to 
the containment $\prec$ defined by $([m],G)\prec([n],H)$ if and only if there is an increasing injection $f:\;[m]\to[n]$ such that $\{f(E)\;|\;E\in G\}=\{f([m])\cap E\;|\;E\in H\}$.
The following conjecture was proposed in \cite{balo_boll_morr_hype}.

\begin{problem}
If $P$ is a hereditary property of ordered hypergraphs, then either $f_P(n)\le c^n$ 
for every $n$ for some constant $c>1$ or one has 
$f_P(n)\ge\sum_{k=0}^{\lfloor n/2\rfloor}{n\choose 2k}k!$ for every $n$. Moreover, the 
lower bound is best possible.
\end{problem}

\noindent
As noted in \cite{balo_boll_morr_hype}, now it is no longer true that in the latter 
case $P$ must contain $\Pi$.

{\bf Tuples of nonnegative integers.} For a fixed $k$ in $\N$, we set $U=\N_0^k$, so 
$U$ contains all $k$-tuples of nonnegative integers. We define the containment by 
$a=(a_1,a_2,\ds,a_k)\prec b=(b_1,b_2,\ds,b_k)$ if and only if $a_i\le b_i$ for every $i$. 
By Higman's theorem, $(U,\prec)$ is wqo.
So there are only countably many downsets. The size 
function $\|\cdot\|$ on $U$ is given by $\|a\|=a_1+a_2+\ds+a_k$. For a downset $P$ 
in $(U,\prec)$, $f_P(n)$ counts all tuples in $P$ whose entries sum up to $n$. 
Stanley \cite{stan_AMM} (see also \cite[Exercise 6 in Chapter 4]{stan1}) 
obtained the following result.

\begin{theorem}[Stanley]\label{Sammthm}
If $P$ is a downset of $k$-tuples of nonnegative integers, then there is a 
rational polynomial $p(x)$ such that $f_P(n)=p(n)$ for large $n$. 
\end{theorem}

\noindent
In fact, the theorem holds for upsets as well because they are complements of 
 downsets and $\#\{a\in\N_0^k\;|\;\|a\|=n\}={n+k-1\choose k-1}$ is a rational 
polynomial in $n$. Jel\'\i nek and Klazar \cite{jeli_klaz} noted that the theorem holds
for the larger class of sets $P\sus\No^k$ that are finite unions of the generalized 
orthants $\{a\in\No^k\;|\;a_i\ge_ib_i, i\in[k]\}$, here $(b_1,\ds,b_k)$ is in $\No^k$ 
and each $\ge_i$ is either $\ge$ or equality $=$; we call such $P$ {\em simple sets}. 
It appears that this generalization 
of Theorem~\ref{Sammthm} provides a unified explanation 
of the exact polynomial results in Theorems~\ref{setpartstir}, \ref{KKthm}, 
\ref{BBWherpro}, and \ref{BBMorgr}, by mapping downsets of structures in a size-preserving 
manner onto simple sets in $\No^k$ for some $k$.

\subsection{Growths of profiles of relational structures}
\label{profily}

In this subsection we mostly follow the survey article of Pouzet \cite{pouz_surv}, 
see also Cameron \cite{came_00}. This approach of relational structures was pioneered by 
Fra\"\i ss\'e \cite{frai_phd, frai_54,frai_rel}. 
A {\em relational structure} $R=(X,(R_i\;|\;i\in I))$ on $X$ is formed 
by the underlying set $X$ and 
relations $R_i\sus X^{m_i}$ on $X$; the sets $X$ and $I$ may be infinite. The {\em size} 
of $R$ is the cardinality of $X$ and $R$ is called finite (infinite) if $X$ is finite 
(infinite). The {\em signature} of $R$ is the list $(m_i\;|\;i\in I)$ 
of arities $m_i\in\N_0$ of the relations $R_i$.  It is {\em bounded} if the numbers 
$m_i$ are bounded and is {\em finite} if $I$ is finite. 
 
Consider two relational structures $R=(X,(R_i\;|\;i\in I))$ and $S=(Y,(S_i\;|\;i\in I))$ 
with the same signature and an injection $f:\;X\to Y$ satisfying, for every $i$ in 
$I$ and every $m_i$-tuple $(a_1,a_2,\ds,a_{m_i})$ in $X^{m_i}$,  
$$
(a_1,a_2,\ds,a_{m_i})\in R_i\iff (f(a_1),f(a_2),\ds,f(a_{m_i}))\in S_i.
$$ 
If such an injection $f$ exists, we say that $R$ is {\em embeddable} in $S$ and 
write $R\prec S$. If in addition $f$ is an identity (in particular, $X\sus Y$), 
$R$ is a {\em substructure} of $S$ and we write $R\prec^*S$. If the injection $f$ is onto $Y$, 
we say that $R$ and $S$ are {\em isomorphic}. 

The {\em age} of a (typically infinite) relational structure $R$ on $X$ is the set $P$ of all 
finite substructures of $R$. Note that the age forms a downset in the poset $(U,\prec^*)$ 
of all finite relational structures with the signature of $R$ whose underlying 
sets are subsets of $X$. The {\em kernel} of $R$ is the set of elements $x$ in $X$ such that 
the deletion of $x$ changes the age. The {\em profile} of $R$ is the unlabeled 
counting function $g_R(n)$ that counts nonisomorphic structures with size $n$ in the 
age of $R$. We get the same function if we replace the age of $R$ by the set $P$ 
of all finite substructures embeddable in $R$ whose underlying sets are $[n]$:
$$
g_R(n)=\#(\{S=([n],(S_i\;|\;i\in I))\;|\;S\prec R\}/\!\sim)
$$ 
where $\sim$ is the isomorphism relation.

The next general result on growth of profiles was obtained by Pouzet \cite{pouz_phd}, 
see also \cite{pouz_surv}. 

\begin{theorem}[Pouzet]\label{pouz}
If $g_R(n)$ is the profile of an infinite relational structure $R$ with bounded signature or 
with finite kernel, then exactly one of the three cases occurs.
\begin{enumerate}
\item The function $g_R(n)$ is constant for large $n$.
\item There are constants $k$ in $\N$ and $0<c<d$ such that $cn^k<g_R(n)<dn^k$ for every $n$.
\item One has $g_R(n)>n^k$ for large $n$ for every constant $k$ in $\N$.
\end{enumerate}
\end{theorem}

\noindent
Case 1 follows from the interesting fact that every infinite $R$ has a nondecreasing profile
(Pouzet \cite {pouz_81}, see \cite[Theorem 4]{pouz_surv} for further discussion) and cases 
2 and 3 were proved in \cite{pouz_phd} (see \cite[Theorems 7 and 42]{pouz_surv}). It is easy
to see (\cite[Theorem 10]{pouz_surv}) that for unbounded signature one can get arbitrarily slowly growing unbounded profiles. Also, it turns out (\cite[Fact 2]{pouz_surv}) that for 
bounded signature and finite-valued profile, one may assume without loss of generality that the signature is finite. 
An infinite graph $G=(\N,E)$ whose components are cliques and that for every $n$ has 
infinitely many components (cliques) of size $n$ shows that for the signature $(2)$ the 
numbers of integer partitions $p_n$ appear as a profile (cf. Theorem~\ref{unlabcount})---for relational structures in general (and unlabeled count) there is no polynomial-exponential 
jump (but cf. Theorem~\ref{BBMtourn}). 

The survey \cite{pouz_surv} contains, besides further results and problems on 
profiles of relational structures, the following attractive conjecture which was partially 
resolved by Pouzet and Thi\'ery \cite{pouz_thie}.
\begin{problem}
In the cases 1 and 2 of Theorem~\ref{pouz} function $g_R(n)$ is a quasipolynomial 
for large $n$. 
\end{problem}

\noindent
As remarked  in \cite{pouz_surv}, since $g_P(n)$ is nondecreasing, if the conjecture holds 
then the leading coefficient in the quasipolynomial must be constant and, in cases 1 and 2, 
$g_R(n)=an^k+O(n^{k-1})$ for some constants $a>0$ and $k$ in $\N_0$ (cf. 
Theorem~\ref{unlabcount} and the following comment). 

Relational structures are quite general in allowing arbitrarily many relations with arbitrary arities and therefore they can accommodate many previously discussed combinatorial 
structures and many more. On the other 
hand, ages of relational structures are less general than downsets of structures, every age
is a downset in the substructure ordering but not vice versa---many downsets of finite 
structures do not come from a single infinite structure (a theorem due to Fra\"\i ss\'e 
\cite[Lemma 7]{pouz_surv}, \cite{came_00} characterizes downsets that are ages). 
An interesting research direction may be to join the general sides of both approaches. 

\section{Five topics in general enumeration}

In this section we review five topics in general enumeration. As we shall see, there are connections
to the results on growth of downsets presented in the previous section. 

\subsection{Counting lattice points in polytopes} 

A {\em polytope} $P$ in $\R^k$ is a convex hull of a finite set of points.
If these points have rational, respectively integral, coordinates, we speak of
{\em rational}, respectively {\em lattice}, polytope. For a polytope $P$ and $n$ in $\N$ 
we consider the dilation $nP=\{nx\;|\;x\in P\}$ of $P$ and the number of lattice points in it, 
$$
f_P(n)=\#(nP\cap\Z^k).
$$
The following useful result was derived by Ehrhart \cite{ehrh} and Macdonald \cite{macd1, macd2}.

\begin{theorem}[Ehrhart, Macdonald]\label{EMthm}
If $P$ is a lattice polytope, respectively rational polytope, and $f_P(n)$ counts lattice points in the dilation $nP$,  
then there is a rational polynomial, respectively rational quasipolynomial, $p(x)$ such that $f_P(n)=p(n)$ 
for every $n$.
\end{theorem}

\noindent
For further refinements and ramifications of this result and its applications see Beck and Robins 
\cite{beck_robi} (also Stanley \cite{stan1}). Barvinok \cite{barv} and Barvinok and Woods 
\cite{barv_wood} developed a beautiful and powerful theory producing polynomial-time algorithms for counting lattice points in rational polytopes. In way of specializations one obtains from it many Wilfian formulas.
We will not say more on it because in its generality it is out of scope of this overview 
(as we said, this not a survey on $\#$P).

\subsection{Context-free languages} 

A {\em language} $P$ is a subset of $A^*$, the infinite set of finite words over a finite alphabet $A$. The natural size function $|\cdot|$ measures 
length of words and 
$$
f_P(n)=\#\{u\in P\;|\;|u|=n\}
$$
is the number of words in $P$ with length $n$. In this subsection the alphabet $A$ is always 
finite, thus $f_P(n)\le |A|^n$.

We review the definition of context-free languages, for further information on
(formal) languages see Salomaa \cite{salom}. A {\em context-free grammar} is a 
quadruple $G=(A,B,c,D)$ where $A,B$ are finite disjoint sets, $c\in B$ (starting 
variable) and $D$ (production rules) is a finite set of pairs $(d,u)$ where $d\in B$ and
$u\in(A\cup B)^*$. A {\em rightmost derivation} of a word $v\in(A\cup B)^*$ in $G$ is a 
sequence of words $v_1=c,v_2,\ds,v_r=v$ in $(A\cup B)^*$ such that $v_i$ is obtained from
$v_{i-1}$ by replacing the rightmost occurrence of a letter $d$ from $B$ in $v_{i-1}$ by 
the word $u$, according to some production rule $(d,u)\in D$. (Note that no $v_i$ with $i<r$
is in $A^*$.) We let $L(G)$ denote the set of words in $A^*$ that have rightmost derivation 
in $G$. If in addition every $v$ in $L(G)$ has a unique rightmost derivation in $G$, then 
$G$ is an {\em unambiguous context-free grammar}.
A language $P\sus A^*$ is {\em context-free} if $P=L(G)$ for a context-free grammar 
$G=(A,B,c,D)$. $P$ is, in addition, {\em unambiguous} if it can be generated 
by an unambiguous context-free grammar. If $P$ is context-free but not unambiguous, we say
that $P$ is {\em inherently ambiguous}. We associate with a context-free grammar $G=(A,B,c,D)$
a digraph $H(G)$ on the vertex set $B$ by putting an arrow $d_1\to d_2$, $d_i\in B$, if and 
only if there is a production rule $(d_1,u)\in D$ such that $d_2$ appears in $u$. We call 
a context-free language {\em ergodic} if it can be generated by a context-free grammar 
$G$ such that the digraph $H(G)$ is strongly connected.

Chomsky and Sch\"utzenberger \cite{chom_schu} obtained the following important result.

\begin{theorem}[Chomsky and Sch\"utzenberger]\label{ChSchthm}
If $P$ is an unambiguous context-free language and $f_P(n)$ counts words of length $n$ in $P$, then the generating function 
$$
F(x)=\sum_{n\ge 0}f_P(n)x^n
$$ 
of $P$ is algebraic over $\Q(x)$.
\end{theorem}

The algebraicity of a power series $F(x)=\sum_{n\ge 0}a_nx^n$ with $a_n$ in $\N_0$ 
has two important practical corollaries for the counting sequence $(a_n)_{n\ge 1}$. 
First, as we already mentioned, it is holonomic. Second, it has a nice asymptotics. 
More precisely, $F(x)$ determines a function analytic in a neighborhood of $0$ and if 
$F(x)$ is not a polynomial, it has a finite radius of convergence $\rho$, $0<\rho\le 1$, 
and finitely many (dominating) singularities on the circle of convergence $|x|=\rho$. 
In the case of single dominating singularity we have
$$
a_n\sim cn^{\alpha}r^n
$$
where $c>0$ is in $\R$, $r=1/\rho\ge 1$ is an algebraic number, and 
the exponent $\alpha$ is in $\Q\backslash\{-1,-2,-3,\ds\}$ (if $F(x)$ is rational, then 
$\alpha$ is in $\N_0$). For example, for Catalan 
numbers $C_n=\frac{1}{n+1}{2n\choose n}$ and their generating function 
$C(x)=\sum_{n\ge 0}C_nx^n$ we have 
$$
xC^2-C+1=0\ \mbox{ and }\ C_n\sim\pi^{-1/2}n^{-3/2}4^n.
$$ 
For more general results on asymptotics of coefficients of algebraic power series 
see Flajolet and Sedgewick \cite[Chapter VII]{flaj_sedg}. 

Flajolet \cite{flaj} used Theorem~\ref{ChSchthm} to prove the inherent ambiguity of certain context-free 
languages. For further information on rational and algebraic power series in enumeration and
their relation to formal languages see Barcucci et al. \cite{barc_al}, Bousquet-M\'elou \cite{mbm_STACS, mbm_ICM}, 
Flajolet and Sedgewick \cite{flaj_sedg} and Salomaa and Soittola \cite{salo_soit}.

How fast do context-free languages grow? Trofimov \cite{trof} proved for them a polynomial to exponential jump.
 
\begin{theorem}[Trofimov]\label{Tthm}
If $P\sus A^*$ is a context-free language over the alphabet $A$, then either 
$f_P(n)\le |A|n^k$ for every $n$
or $f_P(n)>c^n$ for large $n$, where $k>0$ and $c>1$ are constants. 
\end{theorem}

\noindent
Trofimov proved that in the former case in fact $P\sus w_1^*w_2^*\ds w_k^*$ 
for some $k$ words $w_i$ in $A^*$. Later this theorem was independently rediscovered 
by Incitti \cite{inci} and Bridson and Gilman \cite{brid_gilm}.
D'Alessandro, Intrigila and Varricchio \cite{ales_intr_varr}
show that in the former case the function $f_P(n)$ is in fact a quasipolynomial $p(n)$ for 
large $n$ (and that $p(n)$ and the bound on $n$ can be effectively determined from $P$).

Recall that for a language $P\sus A^*$ the growth constant is defined as 
$c(P)=\limsup f_P(n)^{1/n}$. $P$ is {\em growth-sensitive} if $c(P)>1$ and 
$c(P\cap Q)<c(P)$ for every downset $Q$ in $(A^*,\prec)$, where $\prec$ is the subword 
ordering, such that $P\cap Q\ne P$. In other words, forbidding any word $u$ such that 
$u\prec v\in P$ for some $v$ as a subword results in a significant decrease in growth. 
Yet in other words, in growth-sensitive languages an analogue of Marcus--Tardos theorem 
(Theorem~\ref{MTthm}) holds. In a series of papers Ceccherini-Silberstein, Mach\`\i\ and 
Scarabotti \cite{cecc_mach_scar}, Ceccherini-Silberstein and Woess 
\cite{cecc_woes1, cecc_woes2}, and Ceccherini-Silberstein \cite{cecc} the following
theorem on growth-sensitivity was proved. 

\begin{theorem}[Ceccherini-Silberstein, Mach\`\i, Scarabotti, Woess]
Every unambiguous context-free language $P$ that is ergodic and has $c(P)>1$ is 
growth-sensitive.
\end{theorem}

\noindent
See \cite{cecc_woes1} for extensions of the theorem to the ambiguous case and 
\cite{cecc_mach_scar} for the more elementary case of regular languages.

\subsection{Exact counting of regular and other graphs}

We consider finite simple graphs and, for a given set $P\sus\N_0$, the counting function
($\deg(v)=\deg_G(v)$ is the degree of a vertex $v$ in $G$, the number of incident edges)
$$
f_P(n)=\#\{G=([n],E)\;|\;\deg_G(v)\in P\mbox{ for every $v$ in $[n]$}\}.
$$ 
For example, for $P=\{k\}$ we count labeled $k$-regular graphs on $[n]$.
The next general theorem was proved by Gessel \cite[Corollary 11]{gess}, by means of symmetric functions 
in infinitely many variables.

\begin{theorem}[Gessel]\label{Gthm}
If $P$ is a finite subset of $\N_0$ and $f_P(n)$ counts labeled graphs on $[n]$ with 
all degrees in $P$, then the sequence $(f_P(n))_{n\ge 1}$ is holonomic.
\end{theorem}

\noindent
This theorem was conjectured and partially proved for the $k$-regular case, $k\le 4$, 
by Goulden and Jackson \cite{goul_jack}. As remarked in \cite{gess}, the theorem holds 
also for graphs with multiple edges and/or loops.
Domoco\c{s} \cite{domo} extended it to 3-regular and 3-partite hypergraphs (and remarked 
that Gessel's method works also for general $k$-regular and $k$-partite hypergraphs). 
For more information see also Mishna \cite{mish_phd,mish}.

Consequently, the numbers of labeled graphs with degrees in fixed finite set have 
Wilfian formula. Now we demonstrate this directly by a more generally applicable argument. 
For $d$ in $\N_0^{k+1}$, we say that a graph $G$ is a $d$-{\em graph} if 
$|V(G)|=d_0+d_1+\ds+d_k$, $\deg(v)\le k$ 
for every $v$ in $V(G)$ and exactly $d_i$ vertices in $V(G)$ have degree $i$. 
Let
$$
p(d)=p(d_0,d_1,\ds,d_k)=\#\{G=([n],E)\;|\;\mbox{$G$ is a $d$-graph}\}
$$ 
be the number of labeled $d$-graphs with vertices $1,2,\ds,n=d_0+\ds+d_k$.

\begin{proposition}
For fixed $k$, the list of numbers 
$$
(p(d)\;|\;d\in\N_0^{k+1}, d_0+d_1+\ds+d_k=m\le n)
$$ 
can be generated in time polynomial in $n$.
\end{proposition}
\duk
A natural idea is to construct graphs $G$ with $d_i$ vertices of degree $i$ by adding 
vertices $1,2,\ds,n$ one by one, keeping track of the numbers $d_i$.
In the first phase of the algorithm we construct an auxiliary $(n+1)$-partite graph 
$$
H=(V_0\cup V_1\cup\ds\cup V_n,E)
$$ 
where we start with $V_m$ consisting of all $(k+1)$-tuples $d=(d_0,d_1,\ds,d_k)$ in $\N_0^{k+1}$ 
satisfying $d_0+\ds+d_k=m$ and the edges will go only between $V_m$ and $V_{m+1}$. 
An edge joins $d\in V_m$ with $e\in V_{m+1}$ if and only if there exist numbers $\Delta_0,\Delta_1,\ds,\Delta_{k-1}$ 
in $\N_0$ such that: $0\le\Delta_i\le d_i$ for $0\le i\le k-1$, 
$r:=\Delta_0+\ds+\Delta_{k-1}\le k$, and
$$
e_i=d_i+\Delta_{i-1}-\Delta_i\mbox{ for }i\in\{0,1,\ds,k\}\backslash\{r\}
\mbox{ but } e_r=d_r+\Delta_{r-1}-\Delta_r+1
$$  
where we set $\Delta_{-1}=\Delta_k=0$. We omit from $H$ (or better, do not construct 
at all) the vertices $d$ in $V_m$ not reachable from $V_0=\{(0,0,\ds,0)\}$ by a path 
$v_0,v_1,\ds,v_m=d$ with $v_i$ in $V_i$. For example, the $k$ vertices in $V_1$ with 
$d_0=0$ are omitted and only $(1,0,\ds,0)$ remains.
Also, we label the edge $\{e,d\}$ with the $k$-tuple 
$\Delta=(\Delta_0,\ds,\Delta_{k-1})$. (It follows that $\Delta$ and $r$ are uniquely 
determined by $d,e$.) The graph $H$ together with its labels can be constructed in time polynomial in $n$. It 
records the changes of the numbers $d_i$ of vertices with degree $i$ caused by adding to 
$G=([m],E)$ new vertex $m+1$; $\Delta_i$ are the numbers of neighbors of $m+1$
with degree $i$ in $G$ and $r$ is the degree of $m+1$.

In the second phase we evaluate a function $p:\;V_0\cup\ds\cup V_n\to\N$ defined 
on the vertices of $H$ by this inductive rule: $p(0,0,\ds,0)=1$ on $V_0$ 
and, for $e$ in 
$V_m$ with $m>0$,
$$
p(e)=\sum_dp(d)\prod_{i=0}^{k-1}{d_i\choose \Delta_i}
$$
where we sum over all $d$ in $V_{m-1}$ such that $\{d,e\}\in E(H)$ and $\Delta$ is the label of 
the edge $\{d,e\}$. It is easy to see that all values of $p$ can be obtained 
in time polynomial in $n$ and that $p(d)$ for $d$ in $V_m$ is the number of 
labeled $d$-graphs on $[m]$. 
\kduk

\noindent
Now we can in time polynomial in $n$ easily calculate
$$
f_P(n)=\sum_dp(d)
$$
as a sum over all $d=(d_0,\ds,d_k)$ in $V_n$, $k=\max P$, satisfying $d_i=0$ when 
$i\not\in P$. Of course, this algorithm is much less effective than the holonomic recurrence
ensured by (and effectively obtainable by the proof of) Theorem~\ref{Gthm}. 
But by this approach we 
can get Wilfian formula also for some infinite sets of degrees $P$, for example when 
$P$ is an arithmetical progression. (We leave to the reader as a nice exercise to 
count labeled graphs with even degrees.) On the other hand, it seems 
to fail for many classes of graphs, for example, for triangle-free graphs.

\begin{problem}\label{Prtria}
Is there a Wilfian formula for the number of labeled triangle-free graphs on $[n]$? Can 
this number be calculated in time polynomial in $n$?
\end{problem}

\noindent
The problem of enumeration of labeled triangle-free graphs was mentioned by Read in 
\cite[Chapter 2.10]{read_conn}. A quarter century ago,  Wilf \cite{wilf_what} posed the following similar 
problem. 

\begin{problem}
Can one calculate in time polynomial in $n$ the number of unlabeled graphs on $[n]$?
\end{problem}

\subsection{The ubiquitous asymptotics $cn^{-3/2}r^n$}

In many enumerative problems about recursively defined structures, 
for example when counting rooted trees of various kinds, one ends up with asymptotic 
form of the type 
$$
f_P(n)\sim cn^{-3/2}r^n
$$
where $c>0$ and $r>1$ are constants. Bell, Burris and Yeats \cite{bell_burr_yeat} developed 
a remarkable general theory explaining ubiquity of this asymptotic relation. Their results
give practical and general tool for proving this asymptotics by checking certain conditions 
for the operator on power series $\Theta$ that appears in the equation $F=\Theta(F)$ 
expressing the counting problem in terms of the generating function $F=\sum_{n\ge 0}f_P(n)x^n$. These conditions are often easy to check; the reader is referred to \cite{bell_burr_yeat} for 
examples. We state their main result \cite[Theorem 75]{bell_burr_yeat} below as Theorem~\ref{BBYthm}. 
For the complete statement of the theorem it is necessary to introduce several 
notions. Eventually we define two classes ${\cal O}_E$ and ${\cal O}_I$ of operators on power series guaranteeing the asymptotics $cn^{-3/2}r^n$.

By $\R[[x,y]]_{\ge 0}$ we denote the set of bivariate power series with nonnegative 
real coefficients and zero constant term, and similarly for $\R[[x]]_{\ge 0}$ and 
$\Z[[x]]_{\ge 0}$. An {\em operator} is a mapping $\Theta$ from $\R[[x]]_{\ge 0}$ 
to itself. It is {\em integral} if it preserves $\Z[[x]]_{\ge 0}$. We say that $\Theta$ is a 
{\em retro operator} if for every $n$ the coefficient $[x^n]\Theta(F)$ depends only on 
the coefficients $[x^m]F$ with $m<n$, in particular, $[x]\Theta(F)$ is a constant independent 
of $F$. {\em Elementary operator} $\Theta$ is given by a power series $E(x,y)$ in 
$\R[[x,y]]_{\ge 0}$ by $\Theta(F)=E(x,F)$. We also say that $E$ {\em represents} $\Theta$. 
An elementary operator $\Theta=E$ is {\em nonlinear} if $E(x,y)$ is nonlinear in $y$, is {\em bounded} if $[x^n]E(x,x)<c^n$ for every $n$ and a constant $c>1$, 
and is {\em open} if for every $a,b>0$ the convergence $E(a,b)<+\infty$ implies
the convergence $E(a+\varepsilon,b+\varepsilon)<+\infty$ for some $\varepsilon>0$.  
${\cal O}_E$ is the set of bounded and open elementary operators. By ${\cal O}_E^*$ 
we denote the subset of integral operators in ${\cal O}_E$.

The second set of operators ${\cal O}_I$ is defined as a closure, under the operation
of scalar multiplication by positive reals and the binary operations of addition $(+)$, multiplication $(\cdot)$, and composition $(\circ)$ of operators, of the set of 
{\em base operators}
$$
{\cal O}_E^*\cup\{\Theta_{H,M}\;|\;H\in\{S,D,C,I\}, M\sus\N\}
$$
where for $H\in\{D,C\}$ we allow only finite sets $M$ or sets satisfying
$\sum_{m\in M}1/m=+\infty$. The operators $\Theta_{H,M}$ are defined as follows. 
For a permutational group $H$ acting on $[m]$, P\'olya's cycle index polynomial is
$$
Z(H,x_1,\ds,x_m)=\frac{1}{|H|}\sum_{\sigma\in H}\prod_{j=1}^mx_j^{\sigma_j}
$$
where $\sigma_j$ is the number of $j$-cycles in the decomposition of $\sigma$ into disjoint
cycles. In the role of $H$ we take four families of permutational groups: 
$S_m$ (symmetric group of order $m!$), $D_m$ (dihedral group of order $2m$), $C_m$ (cyclic group of order $m$) and $I_m$ (identity group of order $1$). For $m$ in $\N$ we define four operators 
$\Theta_{S,m}$, $\Theta_{D,m}$, $\Theta_{C,m}$, and $\Theta_{I,m}$, by
$$
\Theta_{S,m}(F)=Z(S_m,F(x),F(x^2),\ds,F(x^m)),
$$
and by analogous expressions for the other three operators. For $M\sus\N$ we set
$$
\Theta_{S,M}=\sum_{m\in M}\Theta_{S,m}
$$
and similarly for the cases $D, C, I$.  This completes the definition of the base operators 
$\Theta_{S,M}$, $\Theta_{D,M}$, $\Theta_{C,M}$, and $\Theta_{I,M}$. Note that each of them, 
as well as each operator in ${\cal O}_E^*$, is integral. It follows that every operator in
${\cal O}_I$ is integral.

The base operators $\Theta_{S,M}$, $\Theta_{D,M}$, $\Theta_{C,M}$, and $\Theta_{I,M}$ 
correspond, respectively, to the combinatorial construction of multiset, cycle, directed
cycle and sequence, with $M$ listing the allowed numbers of components. For example, 
the sequence operator $\Theta_{I,M}$ is simply the elementary integral operator
$E(x,y)=\sum_{m\in M}y^m$. See Flajolet and 
Sedgewick \cite{flaj_sedg} for these construction and their relations to generating 
functions, and Bergeron, Labelle and Leroux \cite{species} for a more abstract and general 
approach.

Following the inductive definition of ${\cal O}_I$, one associates with a given $\Theta$ in 
${\cal O}_I$ and $F$ in $\R[[x]]_{\ge 0}$ a canonical power series $E(x,y)=E^{\Theta,F}(x,y)$ 
in $\R[[x,y]]_{\ge 0}$ {\em representing $\Theta$ at $F$}, which means that 
$\Theta(F)=E^{\Theta,F}(x,F)$. We say that $\Theta$ in ${\cal O}_I$ is {\em nonlinear} if 
$E^{\Theta,x}(x,y)$ is nonlinear in $y$. 

We are ready to state the theorem. The lower indices 
in $E_x$ and $E_{yy}$ indicate partial derivatives.

\begin{theorem}[Bell, Burris and Yeats]\label{BBYthm}
Suppose that $\Theta$ is a nonlinear retro operator in ${\cal O}_E$, respectively in 
${\cal O}_I$, and that the power series $A$ in $\R[[x]]_{\ge 0}$, respectively in 
$\Z[[x]]_{\ge 0}$, diverges in its radius of convergence. Then the equation
$$
F=A+\Theta(F)
$$
has a unique solution $F$ in $\R[[x]]_{\ge 0}$, respectively in $\Z[[x]]_{\ge 0}$, and the coefficients of $F$ have the following asymptotic form. There are constants $d\in\N_0$, $q\in\N$, $c>0$, and $r>1$ such that
$[x^n]F=0$ for $n\not\equiv d\;{\rm mod}\;q$ but
$$
[x^n]F\sim cn^{-3/2}r^n
$$
as $n\to\infty$ for $n\equiv d\;{\rm mod}\;q$. Moreover, $r=1/\rho$ where $\rho$ is the 
radius of convergence of $F$ and
$$
c=q\cdot\sqrt{\frac{\rho\cdot E_x(\rho,F(\rho))}{2\pi\cdot E_{yy}(\rho,F(\rho))}}
$$ 
where $E=E(x,y)$ in $\R[[x,y]]_{\ge 0}$ represents the operator 
$L:\;G\mapsto A+\Theta(G)$, respectively $E$ represents $L$ at $F$.
\end{theorem} 

\noindent
An important combinatorial construction left out from the definition of ${\cal O}_I$ is 
the set construction, with corresponding operators $\Theta_{\pm S,M}$ and $\Theta_{\pm S,m}$
where 
$$
\Theta_{\pm S,m}(F)=Z(S_m,F(x),-F(x^2),F(x^3),\ds,(-1)^{m+1}F(x^m)).
$$
Negative signs are troublesome because the nonnegativity of coefficients of power series 
is crucial at several steps of proof of Theorem~\ref{BBYthm}. Besides other interesting problems,
the authors of \cite{bell_burr_yeat} pose the following one 
(\cite[Q2 in Section 8]{bell_burr_yeat}, \cite[Section 6.2]{bell_burr_yeat}).

\begin{problem}
Can the set operator $\Theta_{\pm S,M}$ be adjoined to the base operators in the definition of 
${\cal O}_I$ so that the universal asymptotics $cn^{-3/2}r^n$ still follows?
\end{problem}

\subsection{Ultimate modular periodicity}
\label{logika}

One general aspect of counting functions $f_P(n)$ not touched so far  
is their modular behavior. For given modulus $m$ in $\N$, what can be said about 
the sequence of residues $(f_P(n)\ {\rm mod}\ m)_{n\ge 1}$. Before presenting a rather 
general result in this area, we motivate it by two examples. 

First, Bell numbers $B_n$ counting partitions of $[n]$. Recall that 
$$
\sum_{n\ge 0}B_nx^n=\sum_{k=0}^{\infty}\frac{x^k}{(1-x)(1-2x)\ds(1-kx)}.
$$
Reducing modulo $m$ we get, denoting $v(x)=(1-x)(1-2x)\ds(1-(m-1)x)$, 
\begin{eqnarray*}
\sum_{n\ge 0}B_nx^n&\equiv_m&\sum_{j=0}^{\infty}\frac{x^{mj}}{v(x)^j}
\sum_{i=0}^{m-1}\frac{x^i}{(1-x)(1-2x)\ds(1-ix)}\\
&=&\frac{1}{1-x^m/v(x)}\sum_{i=0}^{m-1}\frac{x^i}{(1-x)(1-2x)\ds(1-ix)}\\
&=&\frac{a(x)}{v(x)-x^m}=\frac{a(x)}{1+b_1x+\ds+b_mx^m}\\
&=&\sum_{n\ge 0}c_nx^n
\end{eqnarray*}
where $a(x)\in\Z[x]$ has degree at most $m-1$ and $b_i$ are integers, $b_m=-1$. 
Thus the sequence of integers $(c_n)_{n\ge 1}$ satisfies for $n>m$ the linear recurrence 
$c_n=-b_1c_{n-1}-\ds-b_mc_{n-m}$ of order $m$. 
By the pigeonhole principle, the sequence $(c_n\ {\rm mod}\ m)_{n\ge 1}$ is periodic for 
large $n$. Since $B_n\equiv c_n\ {\rm mod}\ m$, the sequence $(B_n\ {\rm mod}\ m)_{n\ge 1}$ is periodic for large $n$ as well. (For modular periods of Bell numbers see Lunnon, Pleasants and Stephens \cite{lunn_plea_step}).

Second, Catalan numbers $C_n=\frac{1}{n+1}{2n\choose n}$ counting, for example, noncrossing partitions of $[n]$. The shifted version $D_n=C_{n-1}=\frac{1}{n}{2n-2\choose n-1}$ satisfies 
the recurrence $D_1=1$ and, for $n>1$,
$$
D_n=\sum_{i=1}^{n-1}D_iD_{n-i}=2\sum_{i=1}^{\lfloor n/2\rfloor-1}D_iD_{n-i}+
\sum_{i=\lfloor n/2\rfloor}^{\lceil n/2\rceil}D_iD_{n-i}.
$$
Thus, modulo $2$, $D_1\equiv 1$, $D_n\equiv 0$ for odd $n>1$ and $D_n\equiv D_{n/2}^2
\equiv D_{n/2}$ for even $n$. It follows that $D_n\equiv 1$ if and only if $n=2^m$ and that 
the sequence $(C_n\ {\rm mod}\ 2)_{n\ge 1}$ has $1$'s for $n=2^m-1$ and $0$'s elsewhere. 
In particular, it is not periodic for large $n$. (For modular behavior of Catalan numbers 
see Deutsch and Sagan \cite{deut_saga} and Eu, Liu and Yeh \cite{eu_liu_yeh}).

Bell numbers come out as a special case of a general setting. Consider a {\em relational system}
which is a set $P$ of relational structures $R$ with the same finite signature and underlying 
sets $[n]$ for $n$ ranging in $\N$. We say that $P$ is {\em definable in MSOL}, monadic 
second-order logic, if $P$ coincides with the set of finite models (on sets $[n]$) of a closed 
formula $\phi$ in MSOL. (MSOL has, in addition to the language of the first-order 
logic, variables $S$ for sets of elements, which can be quantified by 
$\forall,\exists$, and atomic formulas of the type $x\in S$; see Ebbinghaus and 
Flum \cite{ebbi_flum}.) Let $f_P(n)$ be the number of relational structures in $P$ on 
the set $[n]$, that is, the number of models of $\phi$ on $[n]$ when $P$ is defined by $\phi$. 
For example, the (first-order) formula $\phi$ given by ($a,b,c$ are variables for
elements, $\sim$ is a binary predicate)
$$
\forall a,b,c: (a\sim a)\;\&\;(a\sim b\Rightarrow b\sim a)\;\&\;((a\sim b\;\&\;b\sim c)
\Rightarrow a\sim c)
$$
has as its models equivalence relations and $f_P(n)=f_{\phi}(n)=B_n$, the Bell numbers.

Let us call a sequence $(s_1,s_2,\ds)$ {\em ultimately periodic} if it is periodic for 
large $n$: there are constants $p,q$ in $\N$ such that $s_{n+p}=s_n$ whenever $n\ge q$.
Specker and Blatter (\cite{blat_spec1,blat_spec2,blat_spec3}, see also Specker \cite{spec})
proved the following remarkable general theorem.

\begin{theorem}[Specker and Blatter]\label{SBthm}
If a relational system $P$ definable in MSOL uses only unary and binary relations and
$f_P(n)$ counts its members on the set $[n]$, then for every $m$ in $\N$ the sequence 
$(f_P(n)\ {\rm mod}\ m)_{n\ge 1}$ is ultimately periodic.
\end{theorem}

\noindent
The general reason for ultimate modular periodicity in Theorem~\ref{SBthm} is the same as 
in our example with $B_n$, residues satisfy a linear recurrence with constant 
coefficients. Fischer \cite{fisc} constructed counterexamples to the theorem for quaternary 
relations, see also Specker \cite{speck05}. Fischer and Makowski \cite{fisc_mako} 
extended Theorem~\ref{SBthm} to CMSOL (monadic second-order logic with modular counting) and to 
relations with higher arities when vertices have bounded degrees. 
Note that regular graphs and triangle-free graphs are first-order definable. 
Thus the counting sequences mentioned in Theorem~\ref{Gthm} and 
in Problem~\ref{Prtria} are ultimately periodic to any modulus. More generally, any hereditary
graph property $P=\av(F)$ with finite base $F$ is first-order definable and similarly for other 
structures. Many hereditary properties with 
infinite bases (and also many sets of graphs or structures which are not hereditary) 
are MSOL-definable; this is the case, for example, for forests ($P=\av(F)$ where $F$ is the 
set of cycles) and for planar graphs (use Kuratowski's theorem). To all of them 
Theorem~\ref{SBthm} applies. On the other hand, as the example with Catalan numbers shows, 
counting of ordered structures is in general out of reach of Theorem~\ref{SBthm}. 

A closely related circle of problems is the determination of spectra of relational 
systems $P$ and more generally of finite models; the {\em spectrum of $P$} is the set 
$$
\{n\in\N\;|\;f_P(n)>0\}
$$
---the set of sizes of members of $P$. In several situation it was 
proved that the spectrum is an ultimately periodic subset of $\N$. See Fischer and Makowski \cite{fisc_mako1} (and the references therein), Shelah \cite{shel} and Shelah and Doron 
\cite{shel_doro}. We conclude with a problem posed in \cite{fisc}.

\begin{problem}   
Does Theorem~\ref{SBthm} hold for relational systems with ternary relations?
\end{problem}   

\bigskip\noindent
{\bf Acknowledgments} My thanks go to the organizers of the conference 
Permutations Patterns 2007 in 
St Andrews---Nik Ru\v{s}kuc, Lynn Hynd, Steve Linton, Vince Vatter, Mikl\'os B\'ona, Einar 
Steingr\'\i msson and Julian West---for a very nice conference which gave me opportunity and incentive to write this overview article. I thank also 
V\'\i t Jel\'\i nek for reading the manuscript and an anonymous referee for making many
improvements and corrections in my style and grammar.


\begin{thebibliography}{50}

\bibitem{albe_atki}
M. H. Albert and M. D. Atkinson, Simple permutations and pattern restricted permutations,
{\em Discrete Math.} 300 (2005), 1--15.

\bibitem{albe_atki_brig}
M. H. Albert, M. D. Atkinson and R. Brignall, Permutation classes of polynomial 
growth, {\em Ann. of Combin.}, to appear. 

\bibitem{albe_al}
M. H. Albert,  M. Elder, A. Rechnitzer, P. Westcott and M. Zabrocki, 
On the Stanley-Wilf limit of 
4231-avoiding permutations and a conjecture of Arratia, {\em Adv. in Appl. Math.} 36  (2006), 96--105.

\bibitem{albe_lint}
M. H. Albert and S. A. Linton, Growing at a perfect speed, {\em Combin. Probab. 
Comput.}, submitted.

\bibitem{albe_lint_rusk}
M. H. Albert, S. Linton and N. Ru\v skuc, The insertion encoding of permutations,
{\em Electr. J. Combin.} 12 (2005), Research Paper 47, 31 pp.

\bibitem{alek}
V. E. Alekseev, Range of values of entropy of hereditary classes of graphs,
{\em Diskret. Mat.} 4 (1992), 148--157 (Russian); {\em Discrete Math. Appl.} 3 (1993), 
191--199 (English translation).

\bibitem{ales_intr_varr}
F. D'Alessandro, B. Intrigila and S. Varricchio, On the structure of the counting function
of sparse context-free languages, {\em Theor. Comput. Sci.} 356 (2006), 104--117.

\bibitem{allo_shal} J.-P. Allouche and J. Shallit, {\em Automatic Sequences}, 
Cambridge University Press, Cambridge, 2003. 

\bibitem{andr}
G. Andrews, {\em The Theory of Partitions}, Addison-Wesley, Reading MA, 1976.

\bibitem{arra}
R. Arratia, On the Stanley-Wilf conjecture for the number of permutations
avoiding a given pattern, {\em Electron. J. Combin.} 6 (1999) Note 1, 4 pp. 

\bibitem{balo_boll}
J. Balogh and B. Bollob\'as, Hereditary properties of words, {\em  Theor. Inform. Appl.} 
39 (2005), 49--65.

\bibitem{balo_boll_morr_hype}
J. Balogh, B. Bollob\'as and R. Morris, Hereditary properties of partitions, 
ordered graphs and ordered hypergraphs, {\em Europ. J. Combin.} 27 (2006), 1263--1281.

\bibitem{balo_boll_morr_orgr}
J. Balogh, B. Bollob\'as and R. Morris, Hereditary properties of ordered graphs, 
in: M. Klazar, J. Kratochv\'\i l, M. Loebl, J. Matou\v sek, R. Thomas and P. Valtr (eds.), 
{\em Topics in Discrete Mathematics (special edition for J. Ne\v{s}et\v{r}il)}, Springer, 
2006, pp. 179--213.

\bibitem{balo_boll_morr_tour}
J. Balogh, B. Bollob\'as and R. Morris, Hereditary properties of tournaments, 
{\em Electr. J. Combin.} 14 (2007), Research Paper 60, 25 pp.

\bibitem{balo_boll_morr_pose}
J. Balogh, B. Bollob\'as and R. Morris, Hereditary properties of combinatorial 
structures: Posets and oriented graphs, {\em J. Graph Theory} 56 (2007), 311--332.

\bibitem{balo_boll_saks_sos}
J. Balogh, B. Bollob\'as, M. Saks and V. T. S\'os, The unlabeled speed of a 
hereditary graph property, submitted.

\bibitem{balo_boll_simo}
J. Balogh, B. Bollob\'as and M. Simonovits, The number of graphs without forbidden 
subgraphs, {\em J. Combin. Theory, Ser. B} 91 (2004), 1--24. 

\bibitem{balo_boll_wein_spee}
J. Balogh, B. Bollob\'as and D. Weinreich, The speed of hereditary properties of 
graphs, {\em J. Combin. Theory, Ser. B} 79 (2000), 131--156.

\bibitem{balo_boll_wein_penu}
J. Balogh, B. Bollob\'as and D. Weinreich, The penultimate range of growth for 
graph properties,  {\em Europ. J. Combin.} 22 (2001), 277--289.

\bibitem{balo_boll_wein_meas}
J. Balogh, B. Bollob\'as and D. Weinreich, Measures on monotone properties of 
graphs, {\em Discrete Appl. Math.} 116 (2002), 17--36. 

\bibitem{balo_boll_wein_bell}
J. Balogh, B. Bollob\'as and D. Weinreich, A jump to the bell number for 
hereditary graph properties, {\em J. Combin. Theory, Ser. B} 95 (2005), 29--48.

\bibitem{barc_al}
E. Barcucci, A. Del Lungo, A. Frossini and S. Rinaldi, From rational functions to
regular languages, in: {\em Proceeding of FPSAC'00}, Springer, 2000, pp. 633--644.

\bibitem{barv}
A. Barvinok, The complexity of generating functions for integer points in polyhedra 
and beyond, in: {\em International Congress of Mathematicians. Vol. III}, 
Eur. Math. Soc., Z\"urich, 2006, pp. 763--787.

\bibitem{barv_wood}
A. Barvinok and K. Woods, Short rational generating functions for lattice point problems,  
{\em J. Amer. Math. Soc.} 16 (2003), 957--979.

\bibitem{beck_robi}
M. Beck and S. Robins, {\em Computing the Continuous Discretely. Integer-point 
Enumeration in Polyhedra}, Springer, 2007. 

\bibitem{bell_burr_yeat}
J. P. Bell, S. N. Burris and K. A. Yeats, Counting rooted trees: the universal law 
$t(n)\sim C\rho\sp {-n} n\sp {-3/2}$, {\em Electron. J. Combin.} 13 (2006), 
Research Paper 63, 64 pp.

\bibitem{species}
F. Bergeron, G. Labelle and P. Leroux, {\em Combinatorial Species and Tree-like Structures},
Cambridge University Press, 1998. 

\bibitem{bern_noy_wels}
O. Bernardi, M. Noy and D. Welsh, On the growth rate of minor-closed classes of 
graphs, arXiv:06710.2995. 

\bibitem{blat_spec1}
C. Blatter and E. Specker, Le nombre de structures finies d'une th\'eorie \`a caract\`ere
fini, {\em Sci. Math. Fonds Nat. Rec. Sci. Bruxelles} (1981), 41--44. 

\bibitem{blat_spec2}
C. Blatter and E. Specker, Modular periodicity of combinatorial sequences, 
{\em Abstracts AMS} 4 (1983), 313. 

\bibitem{blat_spec3}
C. Blatter and E. Specker, Recurrence relations for the number of labelled structures
on a finite set, in: E. B\"orger, G. Hasenjager and D. R\"odding (eds.), {\em In Logic and 
Machines: Decision Problems and Complexity}, Lecture Notes in Computer Science, vol. 171, 
Springer, 1984, pp. 43--61. 

\bibitem{boll_ICM}
B. Bollob\'as, Hereditary properties of graphs: asymptotic enumeration, global 
structure, and colouring, in: {\em Proceedings of the International Congress of 
Mathematicians. Vol. III}, Berlin 1998, {\em Doc. Math. J. DMV} Extra Vol. III
(1998), 333--342.  

\bibitem{boll_bcc07}
B. Bollob\'as,  Hereditary and monotone properties of combinatorial structures, 
in: A. Hilton and J. Talbot (eds.), {\em Surveys in Combinatorics 2007}, 
Cambridge University Press, 2007, pp. 1--40. 

\bibitem{boll_thoma_proj}
B. Bollob\'as and A. Thomason, Projection of bodies and hereditary properties 
of hypergraphs, {\em J. London Math. Soc.} 27 (1995), 417--424.

\bibitem{boll_thoma_here}
B. Bollob\'as and A. Thomason, Hereditary and monotone properties of graphs,
in: R. L. Graham and J. Ne\v{s}et\v{r}il (eds.), {\em The Mathematics of Paul Er\H{o}s 
II}, Springer, 1997, pp. 70--78.

\bibitem{bona_97}
M. B\'ona, Permutations avoiding certain patterns: the case of length $4$ and some 
generalizations, {\em Discrete Math.} 175 (1997), 55--67.

\bibitem{bona_kni}
M. B\'ona, {\em Combinatorics of Permutations}, Chapman $\&$ Hall/CRC, 2004.

\bibitem{boud_pouz}
Y. Boudabbous and M. Pouzet, The morphology of infinite tournaments. Application 
to the growth of their profile, draft, November 2006.  

\bibitem{mbm_STACS}
M. Bousquet-M\'elou, Algebraic generating functions in enumerative combinatorics, 
and context-free languages, in: {\em STACS 2005}, Lecture Notes in Comput. Sci., 
3404, Springer, Berlin, 2005, pp. 18--35.

\bibitem{mbm_ICM}
M. Bousquet-M\'elou, Rational and algebraic series in combinatorial enumeration, in: 
{\em International Congress of Mathematicians. Vol. III}, Eur. Math. Soc., Z\"urich, 
2006, pp. 789--826.

\bibitem{brid_gilm}
M. R. Bridson and R. H. Gilman, Context-free languages of sub-exponential growth,
{\em J. Comput. System Sci.} 64 (2002), 308--310.

\bibitem{brig_hucz_vatt}
R. Brignall, S. Huczynska and V. Vatter, Simple permutations and algebraic generating 
functions, {\em J. Combin. Theory Ser. A} 115 (2008), 423--441. 

\bibitem{brig_rusk_vatt}
R. Brignall, N. Ru\v skuc and V. Vatter, Simple permutations: decidability and unavoidable 
structures, {\em Theor. Comput. Sci.} 391 (2008), 150--163.

\bibitem{burr}
S. N. Burris, {\em Number Theoretic Density and Logical Limit Laws}, AMS, 2001.

\bibitem{came_90}
P. J. Cameron, {\em Oligomorphic Permutation Groups}, Cambridge University Press, 
Cambridge, 1990. 

\bibitem{came_00}
P. J. Cameron, Some counting problems related to permutation groups, {\em Discrete 
Math.} 225 (2000), 77--92.

\bibitem{cecc}
T. Ceccherini-Silberstein, Growth and ergodicity of context-free languages II: the 
linear case, {\em Trans. Amer. Math. Soc.} 359 (2007), 605--618.

\bibitem{cecc_mach_scar} 
T. Ceccherini-Silberstein, A. Mach\`\i\ and F. Scarabotti, On the entropy of regular 
languages, {\em Theor. Comput. Sci.} 307 (2003), 93--102.

\bibitem{cecc_woes1}
T. Ceccherini-Silberstein and W. Woess, Growth and ergodicity of context-free languages,
{\em Trans. Amer. Math. Soc.} 354 (2002), 4597--4625.

\bibitem{cecc_woes2}
T. Ceccherini-Silberstein and W. Woess, Growth-sensitivity of context-free languages,
{\em Theor. Comput. Sci.} 307 (2003), 103--116.

\bibitem{chom_schu}
N. Chomsky and M. P. Sch\"utzenberger, The algebraic theory of context-free languages,
in: {\em Computer Programming and Formal Systems}, North Holland, 1963, 118--161.

\bibitem{comt}
L. Comtet, {\em Advanced Combinatorics}, Reidel, Dordrecht, The Netherlands, 1974.

\bibitem{deut_saga}
E. Deutsch and B. Sagan, Congruences for Catalan and Motzkin numbers and related sequences,
{\em J. Number Theory} 117 (2006), 191--215.

\bibitem{domo}
V. Domoco\c{s}, Minimal coverings of uniform hypergraphs and P-recursiveness, 
{\em Discrete Math.} 159 (1996), 265--271.

\bibitem{ebbi_flum}
H. D. Ebbinghaus and J. Flum, {\em Finite Model Theory}, Springer, 1995.

\bibitem{ehrh}
E. Ehrhart, Sur les poly\`edres rationnels homoth\'etiques \`a $n$ dimensions, 
{\em C. R. Acad. Sci. Paris} 254 (1962), 616--618.

\bibitem{elde_vatt}
M. Elder and V. Vatter, Problems and conjectures presented at the Third International 
Conference on Permutation Patterns (University of Florida, March 7--11, 2005), 
arXiv:math.CO/0505504.

\bibitem{eu_liu_yeh}
S.-P. Eu, S.-Ch. Liu and Y.-N. Yeh, Catalan and Motzkin numbers modulo $4$ and $8$, 
{\em Europ. J. Combin.}, to appear.  

\bibitem{fisc}
E. Fischer, The Specker-Blatter theorem does not hold for quaternary relations,
{\em  J. Combin. Theory Ser. A} 103 (2003), 121--136. 

\bibitem{fisc_mako}
E. Fischer and J. A. Makowski, The Specker-Blatter theorem revisited, in: {\em COCOON'03},
Lecture Notes in Computer Science, vol. 2697, Springer, 2003, pp. 90--101.

\bibitem{fisc_mako1}
E. Fischer and J. A. Makowski, On spectra of sentences of monadic second order logic
with counting, {\em J. Symbolic Logic} 69 (2004), 617--640.

\bibitem{flaj}
P. Flajolet, Analytic models and ambiguity of context-free languages, {\em Theoret. Comput. Sci.}  
49 (1987), 283--309.

\bibitem{flaj_sedg}
P. Flajolet and R. Sedgewick, {\em Analytic Combinatorics}, Cambridge University Press, to appear in 2008.

\bibitem{frai_phd}
R. Fra\"\i ss\'e, {\em Sur quelques classifications des syst\`emes de relations}, Th\`ese, 
Paris, 1953; {\em Alger-Math.} 1 (1954), 35--182.

\bibitem{frai_54}
R. Fra\"\i ss\'e, Sur l'extensions aux relations de quelques propri\'etes des ordres, {\em Ann. 
Sci. Ecole Norm. Sup.} 71 (1954), 361--388.

\bibitem{frai_rel}
R. Fra\"\i ss\'e, {\em Theory of Relations}, North-Holland, Amsterdam, 2000 (second edition).

\bibitem{gess}
I. M. Gessel, Symmetric functions and P-recursiveness, {\em J. Combin. Theory Ser. A} 53 
(1990), 257--285.

\bibitem{goul_jack}
I. P. Goulden and D. M. Jackson, Labelled graphs with small vertex degrees and P-recursiveness, 
{\em SIAM J. Alg. Disc. Meth.} 7 (1986), 60--66. 

\bibitem{grig}
R. Grigorchuk, On the Milnor problem of group growth, {\em Soviet Math. Doklady} 28 (1983),
23--26.

\bibitem{grig_pak}
R. Grigorchuk and I. Pak, Groups of intermediate growth: An introduction for beginners, {\em L'Enseign. Math.}, to appear.

\bibitem{harp}
P. de la Harpe, {\em Topics in Geometric Group Theory}, The University of Chicago Press, 
2000.

\bibitem{higm}
G. Higman, Ordering by divisibility in abstract algebras, {\em Proc. London Math. Soc.}
3 (1952), 326--336. 

\bibitem{hucz_vatt}
S. Huczynska and V. Vatter, Grid classes and the Fibonacci dichotomy for restricted 
permutations, {\em Electron. J. Combin.} 13  (2006), Research Paper 54, 14 pp.

\bibitem{inci}
R. Incitti, The growth function of context-free languages, {\em Theor. Comput. Sci.}
255 (2001), 601--605.

\bibitem{ishi}
Y. Ishigami, The number of hypergraphs and colored hypergraphs with hereditary 
properties, arXiv:0712.0425. 

\bibitem{jeli_klaz} 
V. Jel\'\i nek and M. Klazar, Generalizations of Khovanski\u\i's theorems on 
the growth of sumsets in Abelian semigroups, {\em Adv. Appl. Math.}, to appear.

\bibitem{kais_klaz}
T. Kaiser and M. Klazar, On growth rates of closed sets of hereditary permutation
classes, {\em Electr. J. Combin.} 9 (2002/3), Research Paper 10, 20 pp. 

\bibitem{klaz_stir}
M. Klazar, Counting pattern-free set partitions. I. A generalization of Stirling numbers 
of the second kind, {\em European J. Combin.} 21 (2000), 367--378.

\bibitem{klaz_nch}
M. Klazar, Counting pattern-free set partitions. II. Noncrossing and other 
hypergraphs, {\em Electr. J. Combin.} 7 (2000),  Research Paper 34, 25 pp.

\bibitem{klaz_thcs}
M. Klazar, On the least exponential growth admitting uncountably many closed 
permutation classes, {\em Theor. Comput. Sci.} 321 (2004), 271--281. 

\bibitem{klaz_DS}
M. Klazar, Extremal problems for ordered (hyper) graphs: applications of Davenport-Schinzel sequences, {\em European J. Combin.} 25 (2004), 125--140.

\bibitem{klaz_grrat}
M. Klazar, On growth rates of permutations, set partitions, ordered graphs and other objects, submitted. 

\bibitem{klaz_marc}
M. Klazar and A. Marcus, Extensions of the linear bound in the F\"uredi-Hajnal conjecture,  
{\em Adv. in Appl. Math.} 38 (2007), 258--266.

%\bibitem{kola_prom_roth}
%Ph. G. Kolaitis, H. J. Pr%\"omel and B. L. Rothschild, $K_{l+1}$-free graphs: 
%asymptotic structure and 0--1 laws, {\em Trans. Amer. Math. Soc.} 303 (1987), 
%637--671.

\bibitem{lunn_plea_step}
W. F. Lunnon, P. A. B. Pleasants and N. M. Stephens, Arithmetic properties of Bell numbers 
to a composite modulus. I., {\em  Acta Arith.} 35 (1979), 1--16.

\bibitem{macd1}
I. G. Macdonald, The volume of a lattice polyhedron, {\em Proc. Cambridge Philos. Soc.}
59 (1963), 719--727.

\bibitem{macd2}
I. G. Macdonald, Polynomials associated with finite cell complexes, {\em J. London Math.
Soc.} 4 (1971), 181--192.

\bibitem{macp1}
H. D. Macpherson, Orbits of infinite permutation groups, {\em Proc. London Math. 
Soc. (3)} 51 (1985), 246--284.

\bibitem{macp2}
H. D. Macpherson, Growth rates in infinite graphs and permutation groups, {\em 
Proc. London Math. Soc. (3)} 51 (1985), 285--294.

\bibitem{marc_tard}
A. Marcus and G. Tardos, Excluded permutation matrices and the Stanley--Wilf 
conjecture, {\em J. Combin. Theory, Ser. A} 107 (2004), 153--160.

\bibitem{mari_rado}
D. Marinov and R. Radoi\v ci\'c, Counting 1324-avoiding permutations, {\em Electron. J. Combin.} 
9 (2002/03), Research Paper 13, 9 pp.

\bibitem{mcdi_steg_wels}
C. McDiarmid, A. Steger and D. Welsh, Random graphs from planar and other addable classes, 
in: M. Klazar, J. Kratochv\'\i l, M. Loebl, J. Matou\v sek, R. Thomas and P. Valtr (eds.), 
{\em Topics in Discrete Mathematics (special edition for J. Ne\v{s}et\v{r}il)}, Springer, 
2006, pp. 231--246.

\bibitem{mish_phd}
M. Mishna, {\em Une approche holonome \`a la combinatoire alg\'ebrique}, Doctoral thesis,
Univ. Qu\'ebec \`a Montr\'eal, 2003. 

\bibitem{mish}
M. Mishna, Automatic enumeration of regular objects, {\em J. Integer Seq.} 10 (2007),
Article 07.5.5, 18 pp. 

\bibitem{mors_hedl}
M. Morse and G. A. Hedlund, Symbolic dynamics, {\em Amer. J. Math.} 60 (1938), 
815--866. 

\bibitem{nori_all}
S. Norine, P. Seymour, R. Thomas and P. Wollan, Proper-minor closed families are
small, {\em  J. Combin. Theory, Ser. B} 96 (2006), 754--757.

\bibitem{papa}
Ch. Papadimitriou, {\em Computational Complexity}, Addison-Wesley, 1994.

\bibitem{pouz_phd}
M. Pouzet, {\em Sur la th\'eorie des relations}, Th\'ese d'\'etat, Universit\'e 
Claude-Bernard, Lyon 1, 1978.

\bibitem{pouz_81}
M. Pouzet, Application de la notion de relation presque-encha\i nable au d\'enombrement
des restrictions finies d'une relation, {\em Z. Math. Logik Grundl. Math.} 27 (1981), 
289--332.

\bibitem{pouz_surv}
M. Pouzet, The profile of relations, arXiv:math.CO/0703211.

\bibitem{pouz_thie}
M. Pouzet and N. M. Thi\'ery, Some relational structures with polynomial growth 
and their associated algebras, arXiv:math.CO/0601256.

\bibitem{read_conn}
C. R. Read, Enumeration, in: L. W. Beineke and R. J. Wilson (eds.), 
{\em Graph Connections}, Clarendon Press, Oxford, 1997, pp. 13--33.

\bibitem{minors}
N. Robertson and P. Seymour, Graph minors I--XX, {\em  J. Combin. Theory, Ser. B}, 
1983--2004. 

\bibitem{salom}
A. Salomaa, {\em Formal Languages}, Academic Press, 1973. 

\bibitem{salo_soit}
A. Salomaa and M. Soittola, {\em Automata-theoretic Aspects of Formal Power Series},
Springer, 1978. 

\bibitem{sche_zito}
E. R. Scheinerman and J. Zito, On the size of hereditary classes of graphs, 
{\em J. Combin. Theory, Ser. B} 61 (1994), 16--39. 

\bibitem{shel}
S. Shelah, Spectra of monadic second order sentences, {\em Sci. Math. Japan} 59 (2004),
351--355. 

\bibitem{shel_doro}
S. Shelah and M. Doron, Bounded $m$-ary patch-width are equivalent for $m\ge 3$, 
Shelah's archive, paper no. 865, preprint, 2006. 

\bibitem{simi}
R. Simion, Noncrossing partitions, {\em Discrete Math.} 217 (2000), 367--409.

\bibitem{spec}
E. Specker, Applications of logic and combinatorics to enumeration problems, 
in: E. B\"orger (ed.), {\em Trends in Theoretical Computer Science}, Computer 
Science Press, 1988, pp. 141--169. Reprinted in: {\em Ernst Specker, Selecta},
Birkh\"auser, 1990, pp. 324--350. 

\bibitem{speck05}
E. Specker, Modular counting and substitution of structures, {\em Combin. Probab. 
Comput.} 14 (2005), 203--210.

\bibitem{spen}
J. Spencer, {\em The Strange Logic of Random Graphs}, Springer, 2001.

\bibitem{spie_bona}
D. A. Spielman and M. B\'ona, An infinite antichain of permutations, {\em Electron. 
J. Combin.} 7 (2000), Note 2, 4 pp.

\bibitem{stan_AMM}
R. P. Stanley, Problem E2546, {\em Amer. Math. Monthly} 82 (1975), 756; 
solution {\em  Amer. Math. Monthly} 83 (1976), 813--814. 

\bibitem{stan1}
R. P. Stanley, {\em Enumerative Combinatorics. Vol. 1}, Cambridge University Press, 
Cambridge, 1997 (Corrected reprint of the 1986 original).

\bibitem{stan2}
R. P. Stanley, {\em Enumerative Combinatorics. Vol. 2}, Cambridge University Press, 
Cambridge, 1999.

\bibitem{trof}
V. I. Trofimov, Growth functions of some classes of languages, {\em  Cybernetics} 
17 (1982), 727--731 (translated from {\em Kibernetika} (1981), 9--12).

\bibitem{vatt_schemes}
V. Vatter, Enumeration schemes of restricted permutations, {\em Combin. Probab. Comput.}
17 (2008), 137--159.

\bibitem{vatt_manu}
V. Vatter, Small permutation classes, submitted. 

\bibitem{wilf_what}
H. Wilf, What is an answer?, {\em Amer. Math. Monthly} 89 (1982), 289--292.

\end{thebibliography}
\end{document}